\documentclass[11pt]{article}%
\usepackage{graphicx,subfigure}
\usepackage{amsfonts}
\usepackage{mathrsfs}
\usepackage{amssymb}
\usepackage{chngpage}
\usepackage{amsmath}
\usepackage{graphicx}
\usepackage{amsmath,amssymb}
\usepackage[mathscr]{eucal}
\usepackage{xcolor}%
\providecommand{\U}[1]{\protect\rule{.1in}{.1in}}
\allowdisplaybreaks[4]
\newtheorem{theorem}{Theorem}

\newtheorem{corollary}[theorem]{Corollary}

\newtheorem{definition}[theorem]{Definition}

\newtheorem{lemma}[theorem]{Lemma}

\newtheorem{remark}[theorem]{Remark}

\newenvironment{proof}[1][Proof]{\noindent\textbf{#1.} }{\ \rule{0.5em}{0.5em}}
\numberwithin{equation}{section}
\numberwithin{theorem}{section}

\def\({\left(}
\def\){\right)}

\begin{document}
\centerline{\Large\bf Convergence of solutions }\vspace{0.3cm}
\centerline {\Large\bf  for a reaction-diffusion problem with
fractional Laplacian }


\vskip 16pt\centerline{\bf Jiaohui Xu} \centerline {\footnotesize \it Center
for Nonlinear Studies, School of Mathematics,} \centerline{\footnotesize\it
Northwest University, Xi'an 710127, P. R. China}

\medskip

\centerline{{\bf Tom$\acute{\bf a}$s   Caraballo}\footnote{Corresponding author}
\footnote{ E-mail addresses: jiaxu1@alum.us.es (J. Xu),
caraball@us.es (T. Caraballo), jvalero@umh.es (J. Valero).}}
\centerline{\footnotesize\it Dpto. Ecuaciones Diferenciales y
An$\acute{a}$lisis Num$\acute{e}$rico, } \centerline{\footnotesize\it
Facultad de Matem\'aticas, Universidad de Sevilla,  c/ Tarfia s/n,
41012-Sevilla, Spain}

\medskip

\centerline{\bf Jos\'e Valero}

\centerline{\footnotesize\it  Centro de Investigaci\'on Operativa,
Universidad Miguel Hern\'andez de Elche,} \centerline{\footnotesize\it
Avenida de la Universidad s/n, 03202-Elche, Spain} \medskip

\bigbreak\noindent\textbf{Abstract}{\ } \vskip5pt A kind of nonlocal
reaction-diffusion equations on an unbounded domain containing fractional
Laplacian operator is analyzed. To be precise, we prove the convergence of
solutions of the equation governed by the fractional Laplacian to the
solutions of the classical equation governed by the standard Laplacian, when
the fractional parameter grows to $1$. The existence of global attractors is
investigated as well. The novelty of this paper is concerned with the
convergence of solutions when the fractional parameter varies, which, as far
as the authors are aware, seems to be the first result of this kind of
problems in the literature.

\medskip



\section{Introduction}

For a given initial time $\tau\in\mathbb{R}$, our aim in this paper is to
analyze several interesting aspects related to the following problem on
$\mathbb{R}^{m}$ ($m\geq1$),
\begin{equation}%
\begin{cases}
\dfrac{\partial u}{\partial t}+(-\Delta)^{\gamma}u=f(t,x,u)+h(t,x),\quad
(x,t)\in\mathbb{R}^{m}\mathcal{\times}\left(  \tau,\infty\right)  ,\\
u(\tau,x)=u_{\tau}(x),\quad x\in\mathbb{R}^{m},
\end{cases}
\label{Eq0}%
\end{equation}
where $(-\Delta)^{\gamma}$, $\gamma\in(0,1)$, stands for the fractional
Laplacian operator, $h\in L_{loc}^{2}(\mathbb{R};L^{2}(\mathbb{R}^{m}))$ and
$f$ is a continuous function satisfying some appropriate assumptions as
specified later. Precisely, {under the same assumptions imposed on $f$ and
$h$}, we are interested in studying the convergence of solutions of
\eqref{Eq0}, when $\gamma$ grows to $1^{-}$, to the corresponding ones with
$\gamma=1$, that is, the classical reaction-diffusion problem (see e.g.
\cite{BavinVishik90}, \cite{Wang99}),
\begin{equation}%
\begin{cases}
\dfrac{\partial u}{\partial t}-\Delta u=f(t,x,u)+h(t,x),\quad(x,t)\in
\mathbb{R}^{m}\mathcal{\times}\left(  \tau,\infty\right)  ,\\
u(\tau,x)=u_{\tau}(x),\quad x\in\mathbb{R}^{m}.
\end{cases}
\label{Eq1}%
\end{equation}
The well-posedness of problem \eqref{Eq0} has already been stated (for
instance, in \cite{LuQiWangZhang} it is proved for a stochastic equation which
contains our deterministic model as a particular case), as well as the
existence of global attractors (see also \cite{GalWarma}, \cite[Theorem
2.3]{Wang17} for the case of bounded domains). Several stochastic variants
have also been studied (see \cite{53, 54, 55} and the references
therein), in which the existence and upper-semicontinuity of random attractors
were presented when some parameter in the noise term varies. However, as far
as we are aware, there are no studies yet concerning the convergence of
solutions when the fractional parameter $\gamma$ goes to $1^{-}$. This is an
extremely challenging task which will be tackled in the current paper successfully.

The convergence of $(-\Delta)^{\gamma}u(x)$ to $-\Delta u(x)$, for all
$x\in\mathbb{R}^{m}$ as $\gamma\rightarrow1^{-}$ (see \cite[Proposition
4.4]{Nezza}) when $u\in C_{0}^{\infty}(\mathbb{R}^{m})$ is well known.
Nevertheless, this convergence is far away enough to achieve our goal. In the
end, based on this fact, several novel properties related to fractional
Laplacian operators are established in different phase spaces, which play the
key roles to demonstrate the main result.

The reasons to study this kind of problems in our investigation are as
follows. Fractional problem \eqref{Eq0} and its stationary version have been
considered to illustrate and model the motion of nonlinear deflects in
crystals within the area of dislocation dynamics (see, for example, \cite{B1,
26, 36, T2, T1}). Also, in the phase-field and interfacial dynamics framework,
this equation is usually known as the fractional Allen-Cahn equation (see, for
example, \cite{31,41}). However, it is necessary to emphasize that, to analyze
the problem of anomalous diffusion in physics, probability, finance and other
fields of science, the linear fractional parabolic equation $\partial
_{t}u+(-\Delta)^{\gamma}u=0$ with $\gamma\in(0,1)$ is frequently used instead
of the standard parabolic equation $\partial_{t}u-\Delta u=0$ (see, for
example, \cite{1, WYJ, 33, 40, 42, Xu3, Xu1, Xu2}). There is an interesting
work dealing with different aspects of the normal and anomalous diffusion, see
\cite{47} and the references therein. Therefore, we have a great interest in
knowing if the convergence of the anomalous diffusion problem to the normal
one takes place in a smooth and reasonable way. In our opinion, this would
confirm that the mathematical modeling of the problems is appropriate to
describe the real phenomena.

It is worth highlighting that $\frac{1}{C(m,\gamma)}$ (cf. \eqref{EquivNorm})
is not uniformly bounded with respect to $\gamma\in(0,1)$, see \cite[Corollary
4.2]{Nezza}. Consequently, we cannot use the standard arguments when proving
existence of weak solutions to obtain that $u_{n}\rightarrow u$ strongly in
$L^{2}(\tau,T;L^{2}(\mathbb{R}^{m}))$. The lack of this strong convergence
prevents us from taking the convergence of $f$ even if $f$ is a continuous
function. To overcome this difficulty, we have to pay a price, that is,
suppose that $f$ is a sublinear function so that the solution of the limit
problem is regular enough. In this way, we fill the gap of showing the limit
of $f$ by applying a classical monotone method.

The outline of this paper is as follows. Section \ref{s2} is devoted to recall
the concept of fractional Laplacian operator, set up the problem, introduce
hypotheses and notation, as well as the definition of weak solutions. In
Section \ref{s3}, we perform a rigorous analysis of the fractional Laplacian
operators $(-\Delta)^{\gamma}$ with $\gamma\in(0,1)$, which are the crucial
tools used to prove the convergence of solutions. Section \ref{s4} is fully
dedicated to address the main theorems of the paper, concerning the
convergence of solutions to the equation governed by fractional Laplacian to
the classical reaction-diffusion equations as $\gamma\rightarrow1^{-}$, when
the external term $f$ satisfies a sublinear condition. Finally, Section
\ref{s5} deals with the existence of global attractors to the fractional
reaction-diffusion equations, when $f$ contains more general cases and $h$ is
independent of $t$.

\medskip

\section{Preliminaries}

\label{s2}

In this section, we will recall the concept of fractional Laplacian operator
on $\mathbb{R}^{m}$, enumerate the assumptions imposed on nonlinear terms $f$
and $h$ for the fractional PDEs under investigation and introduce the
definitions of solutions to problems \eqref{Eq0} and \eqref{Eq1}, respectively.

\subsection{Fractional setting}

Let $\mathcal{S}$ be the Schwartz space of rapidly decaying $C^{\infty}$
functions on $\mathbb{R}^{m}$. For any fixed $0<\gamma<1$, the fractional
Laplacian operator $(-\Delta)^{\gamma}$ of $u\in\mathcal{S}$ at the point $x$
is defined by,
\begin{equation}
(-\Delta)^{\gamma}u(x)=-\frac{1}{2}C(m,\gamma)\int_{\mathbb{R}^{m}}%
\frac{u(x+y)+u(x-y)-2u(x)}{|y|^{m+2\gamma}}dy,\qquad x\in\mathbb{R}%
^{m},\label{eq2-1}%
\end{equation}
where $C(m,\gamma)$ is a positive constant given by
\begin{equation}
C(m,\gamma)=\frac{\gamma4^{\gamma}\Gamma(\frac{m+2\gamma}{2})}{\pi^{\frac
{m}{2}}\Gamma(1-\gamma)}.\label{eq2-2}%
\end{equation}
For any $0<\gamma<1$, the fractional Sobolev space $W^{\gamma,2}%
(\mathbb{R}^{m}):=H^{\gamma}(\mathbb{R}^{m})$ is defined by
\[
H^{\gamma}(\mathbb{R}^{m})=\bigg\{u\in L^{2}(\mathbb{R}^{m}):\int%
_{\mathbb{R}^{m}}\int_{\mathbb{R}^{m}}\frac{|u(x)-u(y)|^{2}}{|x-y|^{m+2\gamma
}}dxdy<\infty\bigg\},
\]
endowed with the norm
\[
\Vert u\Vert_{H^{\gamma}(\mathbb{R}^{m})}=\left(  \int_{\mathbb{R}^{m}%
}|u(x)|^{2}dx+\int_{\mathbb{R}^{m}}\int_{\mathbb{R}^{m}}\frac{|u(x)-u(y)|^{2}%
}{|x-y|^{m+2\gamma}}dxdy\right)  ^{\frac{1}{2}}.
\]
From now on, we denote by $|\cdot|$ and $\Vert\cdot\Vert_{p}$ the norms in
$\mathbb{R}^{m}$ and $L^{p}(\mathbb{R}^{m})$ for $p\not =2$, respectively,
whereas we denote the norm and the inner product of $L^{2}(\mathbb{R}^{m})$ by
$\Vert\cdot\Vert$ and $(\cdot,\cdot)$, respectively. With some abuse of
notation, $(\cdot,\cdot)$ will stand also for the pairing between
$L^{q}(\mathbb{R}^{m})$ and $L^{p}(\mathbb{R}^{m})$, where $q$ is the
conjugate number of $p$. Moreover, the Gagliardo semi-norm of $H^{\gamma
}(\mathbb{R}^{m}),$ denoted by $\Vert\cdot\Vert_{\dot{H}^{\gamma}%
(\mathbb{R}^{m})},$ is written as
\begin{equation}
\Vert u\Vert_{\dot{H}^{\gamma}(\mathbb{R}^{m})}^{2}=\int_{\mathbb{R}^{m}}%
\int_{\mathbb{R}^{m}}\frac{|u(x)-u(y)|^{2}}{|x-y|^{m+2\gamma}}dxdy,\quad u\in
H^{\gamma}(\mathbb{R}^{m}).\label{eq23}%
\end{equation}
Thus, $\Vert u\Vert_{H^{\gamma}(\mathbb{R}^{m})}^{2}=\Vert u\Vert^{2}+\Vert
u\Vert_{\dot{H}^{\gamma}(\mathbb{R}^{m})}^{2}$ for all $u\in H^{\gamma
}(\mathbb{R}^{m})$. Note that $H^{\gamma}(\mathbb{R}^{m})$ is a Hilbert space
with the inner product,
\[%
\begin{split}
(u,v)_{H^{\gamma}(\mathbb{R}^{m})} &  =\int_{\mathbb{R}^{m}}u(x)v(x)dx\\
&  ~~+\int_{\mathbb{R}^{m}}\int_{\mathbb{R}^{m}}\frac{(u(x)-u(y))(v(x)-v(y))}%
{|x-y|^{m+2\gamma}}dxdy,\qquad\forall u,v\in H^{\gamma}(\mathbb{R}^{m}).
\end{split}
\]
By \cite{Nezza}, we know that the norm $\Vert u\Vert_{H^{\gamma}%
(\mathbb{R}^{m})}$ is equivalent to $\left(  \Vert u\Vert^{2}+\Vert
(-\Delta)^{\frac{\gamma}{2}}u\Vert^{2}\right)  ^{\frac{1}{2}}$ for $u\in
H^{\gamma}(\mathbb{R}^{m})$. More precisely, we have
\begin{equation}
\Vert u\Vert_{H^{\gamma}(\mathbb{R}^{m})}^{2}=\Vert u\Vert^{2}+\frac
{2}{C(m,\gamma)}\Vert(-\Delta)^{\frac{\gamma}{2}}u\Vert^{2},\qquad\forall u\in
H^{\gamma}(\mathbb{R}^{m}).\label{EquivNorm}%
\end{equation}
We will use the notation $V=H^{1}(\mathbb{R}^{m})$, $H=L^{2}(\mathbb{R}^{m})$,
$V_{\gamma}=H^{\gamma}(\mathbb{R}^{m})$, while their dual spaces are denoted
by $V^{\ast}$, $H^{\ast}$ and $V_{\gamma}^{\ast}$, respectively. Moreover, the
norms in $V$ and $V_{\gamma}$ will be denoted by $\Vert\cdot\Vert_{V}$ and
$\Vert\cdot\Vert_{V_{\gamma}}$. It follows from \cite{Nezza} that the
embeddings $V\subset V_{\gamma_{2}}\subset V_{\gamma_{1}}$ are continuous for
$0<\gamma_{1}\leq\gamma_{2}<1$. Identifying $H$ with its dual, we obtain the
following chain of continuous embeddings
\[
V\subset V_{\gamma}\subset H=H\subset V_{\gamma}^{\ast}\subset V^{\ast}%
,\quad\text{ for }\gamma\in\left(  0,1\right)  .
\]

Let $b_{\gamma}:V_{\gamma}\times V_{\gamma}\rightarrow\mathbb{R}$ be a
bilinear form given by
\begin{equation}
b_{\gamma}(v_{1},v_{2})=\frac{1}{2}C(m,\gamma)\int_{\mathbb{R}^{m}}%
\int_{\mathbb{R}^{m}}\frac{(v_{1}(x)-v_{1}(y))(v_{2}(x)-v_{2}(y))}%
{|x-y|^{m+2\gamma}}dxdy,~~\forall v_{1},v_{2}\in V_{\gamma}, \label{Bilinear}%
\end{equation}
where $C(m,\gamma)$ is the constant in (\ref{eq2-2}). We associate the
operator $A^{\gamma}:V_{\gamma}\rightarrow V_{\gamma}^{\ast}$ with
{$b_{\gamma}$ by
\[
<A^{\gamma}(v_{1}),v_{2}>_{(V_{\gamma}^{\ast},V_{\gamma})}=b_{\gamma}%
(v_{1},v_{2}),\qquad\forall v_{1},v_{2}\in V_{\gamma},
\]
where $<\cdot,\cdot>_{(V_{\gamma}^{\ast},V_{\gamma})}$ is the duality pairing
of $V_{\gamma}^{\ast}$ and $V_{\gamma}$. Similarly, pairing between spaces $V$
and $V^{\ast}$ will be denoted by $<\cdot,\cdot>_{(V^{\ast},V)}$.}

\subsection{Setting of the problem}

Throughout this paper, we assume that $f:\mathbb{R\times R}^{m}\mathbb{\times
R}\rightarrow\mathbb{R}$ is continuous, and such that for all $t$,
$u\in\mathbb{R}$ and $x\in\mathbb{R}^{m}$, $f(t,x,$\textperiodcentered$)\in
C^{1}(\mathbb{R})$ and
\begin{equation}
\frac{\partial f}{\partial u}(t,x,u)\leq\sigma, \label{Lip}%
\end{equation}%
\begin{equation}
f(t,x,u)u\leq-\mu|u|^{2}+\psi_{1}(t,x), \label{Diss}%
\end{equation}%
\begin{equation}
|f(t,x,u)|\leq\psi_{2}(t,x)|u|+\psi_{3}(t,x), \label{Growth}%
\end{equation}
where $\mu>0,\ \sigma\geq0$, $\psi_{1}\in L_{loc}^{1}(\mathbb{R}%
;L^{1}(\mathbb{R}^{m})),\ \psi_{2}\in L_{loc}^{\infty}(\mathbb{R};L^{\infty
}(\mathbb{R}^{m})),\ \psi_{3}\in L_{loc}^{2}(\mathbb{R};L^{2}(\mathbb{R}%
^{m}))$ and $\psi_{i}$ are nonnegative for $i=1,2,3$.

Let us establish the definition of weak solutions to problems (\ref{Eq0}%
)-(\ref{Eq1})

\begin{definition}
Let $\tau\in\mathbb{R}$, $u_{\tau}\in H$ and $\gamma\in(0,1)$. The function
$u\in C([\tau,\infty),H)$ is said to be a weak solution to problem (\ref{Eq0})
if $u(\tau)=u_{\tau}$, $u\in L_{loc}^{2}(\tau,\infty;V_{\gamma})$, $\dfrac
{du}{dt}\in L_{loc}^{2}(\tau,\infty;V_{\gamma}^{\ast})$ and $u$ satisfies, for
every $\xi\in V_{\gamma}$,%
\begin{align}
\frac{d}{dt}\left(  u,\xi\right)   &  +\frac{1}{2}C(m,\gamma)\int%
_{\mathbb{R}^{m}}\int_{\mathbb{R}^{m}}\frac{(u(t,x)-u(t,y))(\xi(x)-\xi
(y))}{|x-y|^{m+2\gamma}}dxdy\nonumber\\
&  =\int_{\mathbb{R}^{m}}\left(  f(t,x,u(t,x))+h(t,x)\right)  \xi\left(
x\right)  dx, \label{EqSolutionGamma}%
\end{align}
in the sense of scalar distributions on $\left(  \tau,\infty\right)  .$
\end{definition}

It is known (see \cite[Theorem 3.2]{LuQiWangZhang}) that problem (\ref{Eq0})
possesses a unique weak solution for each $u_{\tau}\in H$ and $\gamma
\in\left(  0,1\right)  $, which is continuous with respect to the initial
datum $u_{\tau}$. Moreover, it satisfies the energy equation,%
\[
\frac{d}{dt}\left\Vert u\right\Vert ^{2}+C(m,\gamma)\Vert u\Vert_{\dot
{H}^{\gamma}(\mathbb{R}^{m})}^{2}=2\int_{\mathbb{R}^{m}}\left(
f(t,x,u(t,x))+h(t,x)\right)  u\left(  t,x\right)  dx,\quad\text{ for a.a.
}t\geq\tau.
\]

\begin{definition}
Let $\tau\in\mathbb{R}$, $u_{\tau}\in H$ and $\gamma=1$. The function $u\in
C([\tau,\infty),H)$ is said to be a weak solution to problem (\ref{Eq1}) if
$u(\tau)=u_{\tau}$, $u\in L_{loc}^{2}(\tau,\infty;V)$, $\dfrac{du}{dt}\in
L_{loc}^{2}(\tau,\infty;V^{\ast})$ and $u$ satisfies, for every $\xi\in V$,%
\begin{equation}
\frac{d}{dt}\left(  u,\xi\right)  +\int_{\mathbb{R}^{m}}\nabla
u\text{\textperiodcentered}\nabla\xi dx=\int_{\mathbb{R}^{m}}\left(
f(t,x,u(t,x))+h(t,x)\right)  \xi\left(  x\right)  dx, \label{EqSolution}%
\end{equation}
in the sense of scalar distributions on $\left(  \tau,\infty\right)  .$
\end{definition}
As before, it is also well known (see e.g. \cite{Wang2009}) that problem
(\ref{Eq1}) possesses a unique weak solution for every $u_{\tau}\in H$, which
is continuous with respect to the initial datum $u_{\tau}$. Also, it satisfies
the energy equation,
\[
\frac{d}{dt}\left\Vert u\right\Vert ^{2}+2\Vert u\Vert_{V}^{2}=2\int%
_{\mathbb{R}^{m}}\left(  f(t,x,u(t,x))+h(t,x)\right)  u\left(  t,x\right)
dx,\quad\text{ for a.a. }t\geq\tau.
\]

Let us define the function $\widetilde{f}:\mathbb{R\times R}^{m}\mathbb{\times
R}\rightarrow\mathbb{R}$ by
\[
\widetilde{f}(t,x,u)=f(t,x,u)-\sigma u.
\]
It is clear that
\begin{equation}
\frac{\partial\widetilde{f}}{\partial u}(t,x,u)\leq0,\label{Lip2}%
\end{equation}%
\begin{equation}
\widetilde{f}(t,x,u)u\leq-\mu|u|^{2}+\psi_{1}(t,x),\label{Diss2}%
\end{equation}%
\begin{equation}
|\widetilde{f}(t,x,u)|\leq(\psi_{2}(t,x)+\sigma)|u|+\psi_{3}%
(t,x).\label{Growth2}%
\end{equation}
In this way, problems (\ref{Eq0}) and (\ref{eq1}) can be  rewritten as%
\begin{equation}%
\begin{cases}
\dfrac{\partial u}{\partial t}+(-\Delta)^{\gamma}u-\sigma u=\widetilde{f}%
(t,x,u)+h(t,x),\quad(x,t)\in\mathbb{R}^{m}\mathcal{\times}\left(  \tau
,\infty\right)  ,\\
u(\tau,x)=u_{\tau}(x),\quad x\in\mathbb{R}^{m},
\end{cases}
\label{Eq0B}%
\end{equation}%
and
\begin{equation}%
\begin{cases}
\dfrac{\partial u}{\partial t}-\Delta u-\sigma u=\widetilde{f}%
(t,x,u)+h(t,x),\quad(x,t)\in\mathbb{R}^{m}\mathcal{\times}\left(  \tau
,\infty\right)  ,\\
u(\tau,x)=u_{\tau}(x),\quad x\in\mathbb{R}^{m},
\end{cases}
\label{Eq1B}%
\end{equation}
respectively.
Further, we define the $\overline{f}:\mathbb{R\times R}^{m}\mathbb{\times
R}\rightarrow\mathbb{R}$ by
\[
\overline{f}(t,x,v)=\widetilde{f}(t,x,ve^{\sigma t})e^{-\sigma t}.
\]
Consequently, 
\begin{equation}
\frac{\partial\overline{f}}{\partial u}(t,x,v)\leq0,\label{Lip3}%
\end{equation}%
\begin{align}
\overline{f}(t,x,v)v &  =\widetilde{f}(t,x,ve^{\sigma t})ve^{-\sigma
t}\nonumber\\
&  \leq-\mu|v|^{2}+\psi_{1}(t,x)e^{-2\sigma t},\label{Diss3}%
\end{align}%
\begin{equation}
|\overline{f}(t,x,v)|\leq(\psi_{2}(t,x)+\sigma)|v|+\psi_{3}(t,x)e^{-\sigma
t}.\label{Growth3}%
\end{equation}
By performing the change of variable $v=ue^{-\sigma t}$, problems (\ref{Eq0}) and (\ref{eq1})
formally become%
\begin{equation}%
\begin{cases}
\dfrac{\partial v}{\partial t}+(-\Delta)^{\gamma}v=\overline{f}%
(t,x,v)+e^{-\sigma t}h(t,x),\quad(x,t)\in\mathbb{R}^{m}\mathcal{\times}\left(
\tau,\infty\right)  ,\\
v(\tau,x)=e^{-\sigma\tau}u_{\tau}(x),\quad x\in\mathbb{R}^{m},
\end{cases}
\label{Eq0C}%
\end{equation}
and
\begin{equation}%
\begin{cases}
\dfrac{\partial v}{\partial t}-\Delta v=\overline{f}(t,x,v)+e^{-\sigma
t}h(t,x),\quad(x,t)\in\mathbb{R}^{m}\mathcal{\times}\left(  \tau
,\infty\right)  ,\\
v(\tau,x)=e^{-\sigma\tau}u_{\tau}(x),\quad x\in\mathbb{R}^{m}.
\end{cases}
\label{Eq1C}%
\end{equation}

\begin{lemma}
\label{EquivSol}The function $u\left(  t\right)  $ is a weak solution to
problem (\ref{Eq0}) if and only if $v\left(  t\right)  =u(t)e^{-\sigma t}$ is
a weak solution to problem (\ref{Eq0C}).
\end{lemma}

\begin{proof}
Let $u\left(  t\right)  $ be a weak solution to problem (\ref{Eq0}). Then it
is clear that $v\in C([\tau,\infty),H)$, $v(\tau)=e^{-\sigma\tau}u_{\tau}$,
$v\in L_{loc}^{2}(\tau,\infty;V_{\gamma})$, $\dfrac{dv}{dt}\in L_{loc}%
^{2}(\tau,\infty;V_{\gamma}^{\ast})$ and%
\[
\dfrac{dv}{dt}=\dfrac{du}{dt}e^{-\sigma t}-\sigma e^{-\sigma t}u.
\]
Hence,%
\begin{align*}
&  \frac{d}{dt}\left(  v,\xi\right)  +\frac{1}{2}C(m,\gamma)\int%
_{\mathbb{R}^{m}}\int_{\mathbb{R}^{m}}\frac{(v(t,x)-v(t,y))(\xi(x)-\xi
(y))}{|x-y|^{m+2\gamma}}dxdy\\
&  -\int_{\mathbb{R}^{m}}\left(  \overline{f}(t,x,v(t,x))+e^{-\sigma
t}h(t,x)\right)  \xi\left(  x\right)  dx\\
&  =e^{-\sigma t}\frac{d}{dt}\left(  u,\xi\right)  +e^{-\sigma t}\frac{1}%
{2}C(m,\gamma)\int_{\mathbb{R}^{m}}\int_{\mathbb{R}^{m}}\frac
{(u(t,x)-u(t,y))(\xi(x)-\xi(y))}{|x-y|^{m+2\gamma}}dxdy\\
& \quad -e^{-\sigma t}\int_{\mathbb{R}^{m}}\left(
\underset{=f(t,x,u(t,x))}{\underbrace{\sigma u(t,x)+\widetilde{f}%
(t,x,e^{\sigma t}v(t,x))}}+h(t,x)\right)  \xi\left(  x\right)  dx\\
&  =0.
\end{align*}

The converse is proved in a similar way.
\end{proof}

\bigskip

Likewise we prove the analogous statement for problem (\ref{Eq1}).

\begin{lemma}
\label{EquivSol2}The function $u\left(  t\right)  $ is a weak solution to
problem (\ref{Eq1}) if and only if $v\left(  t\right)  =u(t)e^{-\sigma t}$ is
a weak solution to problem (\ref{Eq1C}).
\end{lemma}

\section{Some properties of the fractional Laplacian operator $(-\Delta
)^{\gamma}$}

\label{s3}

We establish in this section several properties for the fractional Laplacian
operator $(-\Delta)^{\gamma}$, which are the essential tools to prove the key
result of this manuscript. Let $C_{0}^{\infty}(X)$ be the space of infinitely
differentiable functions $u:X\rightarrow\mathbb{R}$ with compact support.
Denote by $\mathcal{D}^{\prime}$ the space of distributions of $\mathcal{D}%
=C_{0}^{\infty}((\tau,\tau+T)\times\mathbb{R}^{m})$ and by $<\cdot
,\cdot>_{(\mathcal{D}^{\prime},\mathcal{D})}$ the duality pairing of
$\mathcal{D}^{\prime}$ and $\mathcal{D}$. {Besides, let us denote by $B_{R}$
the ball in $\mathbb{R}^{m}$ centered at $0$ with radius $R.$}

Initially, recall the well-known limit from \cite[Proposition 4.4]{Nezza} that
if $u\in C_{0}^{\infty}(\mathbb{R}^{m}),$ then%
\begin{equation}
\lim_{\gamma\rightarrow1^{-}}(-\Delta)^{\gamma}u(x)=-\Delta u(x),~~~\text{ for
all }x\in\mathbb{R}^{m}\text{.} \label{ConvLapl}%
\end{equation}
Now, we will extend this convergence in several phase spaces.

\begin{lemma}
\label{ConvergLaplacianLp} For any $\tau\in\mathbb{R}$, $T>0$ and $u\in
C_{0}^{\infty}((\tau,\tau+T)\times\mathbb{R}^{m})$, the following statement
holds:
\[
\lim_{\gamma\rightarrow1^{-}}(-\Delta)^{\gamma}u=-\Delta u,~\text{strongly
in}~L^{p}((\tau,\tau+T)\times\mathbb{R}^{m}),~~\forall p\geq1.
\]
In particular, $\lim_{\gamma\rightarrow1^{-}}(-\Delta)^{\gamma}u=-\Delta u$ in
the sense of distributions in $\mathcal{D}^{\prime}.$
\end{lemma}

\begin{proof}
Due to the singularity of the kernel (cf. \eqref{eq2-1}) of the fractional
Laplacian operator $(-\Delta)^{\gamma}$, we will split the estimates into two parts.

\textbf{(i)} On the one hand, by \cite[(4.14)]{Nezza}, we have%
\begin{align}
&  ~~\quad C(m,\gamma)\left\vert \int_{B_{1}}\frac{u(t,x+y)+u(t,x-y)-2u(t,x)}%
{|y|^{m+2\gamma}}dy\right\vert \nonumber\\
&  \leq C(m,\gamma)\frac{\omega_{m-1}\left\Vert u\right\Vert _{C^{2}%
((\tau,\tau+T)\times\mathbb{R}^{m})}}{2(1-\gamma)}\leq R_{1},~\text{ }\forall
t\in(\tau,\tau+T),\ x\in\mathbb{R}^{m},\text{ }\gamma\geq\gamma_{0},
\label{EstB1}%
\end{align}
for some constant $R_{1}>0$ and $\gamma_{0}\in\left(  0,1\right)  $, where
$\omega_{m-1}$ is the $(m-1)$-dimensional measure of the unit sphere $S^{m-1}%
$. In order to obtain the last inequality, we have used the following result
from \cite[Corollary 4.2]{Nezza}:%
\[
\lim_{\gamma\rightarrow1^{-}}\frac{C(m,\gamma)}{1-\gamma}=\frac{4m}%
{\omega_{m-1}}.
\]
On the other hand, since $u\in C_{0}^{\infty}((\tau,\tau+T)\times
\mathbb{R}^{m})$, there exists $R>0$ such that,
\begin{equation}
\int_{B_{1}}\frac{u(t,x+y)+u(t,x-y)-2u(t,x)}{|y|^{m+2\gamma}}dy=0,\quad
\mbox{ for }~x\not \in B_{R}. \label{eq1}%
\end{equation}
Combining \eqref{EstB1} and \eqref{eq1} to define
\[
v_{1}(t,x)=%
\begin{cases}%
\begin{split}
R_{1},\quad &  \mbox{if}~~|x|\leq R,\\
0,~\quad &  \mbox{if}~~|x|>R,\\
&
\end{split}
\end{cases}
\]
it is immediate to know that
\[
C(m,\gamma)\left\vert \int_{B_{1}}\frac{u(t,x+y)+u(t,x-y)-2u(t,x)}%
{|y|^{m+2\gamma}}dy\right\vert \leq v_{1}(t,x),\quad\text{ for a.a. }%
x\in\mathbb{R}^{m},
\]
and $v_{1}\in L^{p}((\tau,\tau+T)\times\mathbb{R}^{m})$.

\textbf{(ii)} Secondly, it follows from the H\"{o}lder inequality that there
exists a positive constant $R_{2}$ such that for all $\gamma\geq\gamma_{0}$
with $\gamma_{0}\in(0,1)$ we have
\begin{align}
&  ~~\quad C(m,\gamma)\left\vert \int_{\mathbb{R}^{m}\backslash B_{1}}%
\frac{u(t,x+y)+u(t,x-y)-2u(t,x)}{|y|^{m+2\gamma}}dy\right\vert \nonumber\\
&  \leq C(m,\gamma)\left(  \int_{\mathbb{R}^{m}\backslash B_{1}}\frac
{1}{|y|^{m+2\gamma}}dy\right)  ^{\frac{p-1}{p}}\left(  \int_{\mathbb{R}%
^{m}\backslash B_{1}}\frac{\left(  u(t,x+y)+u(t,x-y)-2u(t,x)\right)  ^{p}%
}{|y|^{m+2\gamma}}dy\right)  ^{\frac{1}{p}}\nonumber\\
&  \leq R_{2}\left(  \frac{\omega_{m-1}}{2\gamma_{0}}\right)  ^{\frac{p-1}{p}%
}\left(  \int_{\mathbb{R}^{m}\backslash B_{1}}\frac{\left(
u(t,x+y)+u(t,x-y)-2u(t,x)\right)  ^{p}}{|y|^{m+2\gamma_{0}}}dy\right)
^{\frac{1}{p}}:=v_{2}(t,x),\label{EstB1out}%
\end{align}
for all $t\in(\tau,\tau+T)$ and $x\in\mathbb{R}^{m}$. Now, we will check
$v_{2}\in L^{p}((\tau,\tau+T)\times\mathbb{R}^{m})$. In light of the
H\"{o}lder inequality, we deduce that there exist positive constants $R_{3}$
and $R_{4}$ such that,
\[%
\begin{split}
&  \qquad\int_{\tau}^{\tau+T}\int_{\mathbb{R}^{m}}|v_{2}\left(  t,x\right)
|^{p}dxdt\\
&  \leq R_{3}\int_{\mathbb{R}^{m}\backslash B_{1}}\frac{1}{|y|^{m+2\gamma_{0}%
}}\int_{\tau}^{\tau+T}\int_{\mathbb{R}^{m}}\left(
u(t,x+y)+u(t,x-y)-2u(t,x)\right)  ^{p}dxdtdy\\[1ex]
&  \leq R_{4}\left\Vert u\right\Vert _{L^{p}((\tau,\tau+T)\times\mathbb{R}%
^{m})}^{p}\frac{\omega_{m-1}}{2\gamma_{0}}.
\end{split}
\]
Hence, let $v(t,x)=\max\{v_{1}(t,x),v_{2}(t,x)\}$ for a.a. $(t,x)\in(\tau
,\tau+T)\times\mathbb{R}^{m}$, we deduce from \textbf{(i)} and \textbf{(ii)}
that $v\in L^{p}((\tau,\tau+T)\times\mathbb{R}^{m})$, and
\[
\left\vert (-\Delta)^{\gamma}u(t,x)\right\vert \leq v(t,x),\qquad\forall
\gamma\geq\gamma_{0}.
\]
The first statement follows from (\ref{ConvLapl}) and Lebesgue's theorem immediately.

Furthermore, for arbitrary $\varphi\in\mathcal{D}$, we obtain%
\[%
\begin{split}
\left\langle (-\Delta)^{\gamma}u+\Delta u,\varphi\right\rangle _{(\mathcal{D}%
^{\prime},\mathcal{D})}=\int_{\tau}^{\tau+T}\int_{\mathbb{R}^{m}}\left(
(-\Delta)^{\gamma}u(t,x)+\Delta u(t,x)\right)  \varphi(t,x)dxdt
\overset{\gamma\rightarrow1^{-}}{\longrightarrow}0.
\end{split}
\]
The proof of this lemma is complete.
\end{proof}


\bigskip

Now, in a similar way we prove the following result.

\begin{lemma}
\label{ConvergLaplacianLp2}For any $u\in C_{0}^{\infty}(\mathbb{R}^{m})$,
$\lim_{\gamma\rightarrow1^{-}}(-\Delta)^{\gamma}u=-\Delta u$ strongly in
$L^{p}(\mathbb{R}^{m})$ for any $p\geq1$.
\end{lemma}

We next establish the continuity of the operator $\left(  -\Delta\right)
^{\gamma}$ with respect to the parameter $\gamma$.

\begin{lemma}
\label{ConvergLaplacianGamma} For any $\tau\in\mathbb{R}$, $T>0$, $u\in
C_{0}^{\infty}((\tau,\tau+T)\times\mathbb{R}^{m})$ and $\gamma_{0}\in(0,1)$,
the following statement holds:
\[
\lim_{\gamma\rightarrow\gamma_{0}}(-\Delta)^{\gamma}u=\left(  -\Delta\right)
^{\gamma_{0}}u,~\text{strongly in}~L^{p}((\tau,\tau+T)\times\mathbb{R}%
^{m}),~~\forall p\geq1.
\]

\end{lemma}

\begin{proof}
Let $\gamma\rightarrow\gamma_{0}\in\left(  0,1\right)  $. We first show
$\lim_{\gamma\rightarrow\gamma_{0}}(-\Delta)^{\gamma}u(t,x)=(-\Delta
)^{\gamma_{0}}u(t,x)$ for any $(t,x)\in(\tau,\tau+T)\times\mathbb{R}^{m}$. It
is obvious that,%
\[
\frac{u\left(  t,x+y\right)  +u(t,x-y)-2u(t,x)}{\left\vert y\right\vert
^{m+2\gamma}}\overset{\gamma\rightarrow\gamma_{0}}{\longrightarrow}%
\frac{u\left(  t,x+y\right)  +u(t,x-y)-2u(t,x)}{\left\vert y\right\vert
^{m+2\gamma_{0}}},~~\text{ for a.a. }y\in\mathbb{R}^{m}.
\]
For $\gamma\leq\gamma_{0}+\alpha$, we have%
\[
\left\vert \frac{u\left(  t,x+y\right)  +u(t,x-y)-2u(t,x)}{\left\vert
y\right\vert ^{m+2\gamma}}\right\vert \leq\frac{\left\Vert D^{2}u\right\Vert
_{L^{\infty}((\tau,\tau+T)\times\mathbb{R}^{m})}}{\left\vert y\right\vert
^{m+2\left(  \gamma_{0}-1+\alpha\right)  }},~\text{ if }\left\vert
y\right\vert <1,
\]
and for $\gamma\geq\gamma_{0}-\alpha$,
\[
\left\vert \frac{u\left(  t,x+y\right)  +u(t,x-y)-2u(t,x)}{\left\vert
y\right\vert ^{m+2\gamma}}\right\vert \leq\frac{3\left\Vert u\right\Vert
_{L^{\infty}((\tau,\tau+T)\times\mathbb{R}^{m})}}{\left\vert y\right\vert
^{m+2\left(  \gamma_{0}-\alpha\right)  }},~\text{ if }\left\vert y\right\vert
\geq1,
\]
where $D^{2}u$ stands for the matrix of second derivatives of $u.$ Hence,
(\ref{eq2-1}), the continuity of Gamma function for positive part and
Lebesgue's theorem yield that
\[%
\begin{split}
\left(  -\Delta\right)  ^{\gamma}u\left(  t,x\right)   &  =-\frac{1}%
{2}C(m,\gamma)\int_{\mathbb{R}^{m}}\frac{u\left(  t,x+y\right)
+u(t,x-y)-2u(t,x)}{\left\vert y\right\vert ^{m+2\gamma}}dy\\
&  \overset{\gamma\rightarrow\gamma_{0}}{\longrightarrow}-\frac{1}%
{2}C(m,\gamma_{0})\int_{\mathbb{R}^{m}}\frac{u\left(  t,x+y\right)
+u(t,x-y)-2u(t,x)}{\left\vert y\right\vert ^{m+2\gamma_{0}}}dy\\[1ex]
&  =\left(  -\Delta\right)  ^{\gamma_{0}}u\left(  t,x\right)  .
\end{split}
\]

In what follows, we check the above convergence is true in $L^{p}((\tau
,\tau+T)\times\mathbb{R}^{m})$. Take a small $\alpha>0$ such that
$|\gamma-\gamma_{0}|<\alpha$, namely, $\gamma_{0}-\alpha<\gamma<\gamma
_{0}+\alpha$. Let us proceed likewise as in the proof of Lemma
\ref{ConvergLaplacianLp}. Notice that
\begin{equation}%
\begin{split}
&  ~\quad C(m,\gamma)\left\vert \int_{\left\vert y\right\vert <1}%
\frac{u(t,x+y)+u(t,x-y)-2u(t,x)}{|y|^{m+2\gamma}}dy\right\vert \\
&  =C(m,\gamma)\left\vert \int_{\left\vert y\right\vert <1}\frac
{1}{|y|^{m+2\gamma}}\int_{0}^{1}\int_{-r}^{r}D^{2}%
u(t,x+sy)y\text{\textperiodcentered}ydsdrdy\right\vert \\
&  \leq K_{1}\int_{\left\vert y\right\vert <1}\frac{1}{\left\vert y\right\vert
^{m+2\left(  \gamma_{0}-1+\alpha\right)  }}\int_{-1}^{1}\left\vert
D^{2}u(t,x+sy)\right\vert dsdy:=v_{3}(t,x),
\end{split}
\label{eq35}%
\end{equation}
where $K_{1}$ is a positive constant. Thanks to $u\in C_{0}^{\infty}%
((\tau,\tau+T)\times\mathbb{R}^{m})$, there exists a sufficiently large $R$
such that, for any $y\in B_{1}$, $s\in\lbrack-1,1]$ and $t\in\lbrack\tau
,\tau+T]$,
\[
D^{2}u(t,x+sy)=0,\quad\forall x\not \in B_{R}.
\]
Also, we find some positive constant $K_{2}$ satisfying $\left\vert
D^{2}u(t,x+sy)\right\vert \leq K_{2}$. Thus, $v_{3}\in L^{p}((\tau
,\tau+T)\times\mathbb{R}^{m})$.

Similarly, for $y\in\mathbb{R}^{m}\backslash B_{1}$, there exists a positive
constant $K_{3}$ such that
\begin{equation}%
\begin{split}
&  \quad~C(m,\gamma)\left\vert \int_{\left\vert y\right\vert \geq1}%
\frac{u(t,x+y)+u(t,x-y)-2u(t,x)}{|y|^{m+2\gamma}}dy\right\vert \\
&  \leq K_{3}\left\vert \int_{\left\vert y\right\vert \geq1}\frac
{u(t,x+y)+u(t,x-y)-2u(t,x)}{|y|^{m+2\gamma_{0}-2\alpha}}dy\right\vert
:=v_{4}(t,x),
\end{split}
\label{eq35B}%
\end{equation}
where the last inequality holds since
\[
\frac{1}{\left\vert y\right\vert ^{m+2\gamma}}\leq\frac{{1}}{\left\vert
y\right\vert ^{m+2\left(  \gamma_{0}-\alpha\right)  }},~~\text{ if }\left\vert
y\right\vert \geq1.
\]
The function $v_{4}$ belongs to $L^{p}((\tau,\tau+T)\times\mathbb{R}^{m})$.
Indeed, by the H\"{o}lder inequality, we derive
\begin{equation}%
\begin{split}
&  ~\quad\int_{\tau}^{\tau+T}\int_{\mathbb{R}^{m}}v_{4}^{p}(t,x)dxdt\\[1ex]
&  \leq\int_{\tau}^{\tau+T}\int_{\mathbb{R}^{m}}K_{3}^{p}\left(
\int_{\left\vert y\right\vert \geq1}\frac{u(t,x+y)+u(t,x-y)-2u(t,x)}%
{|y|^{m+2\gamma_{0}-2\alpha}}dy\right)  ^{p}dxdt\\
&  \leq K_{3}^{p}\left(  \int_{\left\vert y\right\vert \geq1}\frac
{1}{|y|^{m+2\gamma_{0}-2\alpha}}dy\right)  ^{p-1}\\
&  ~~\times\int_{\left\vert y\right\vert \geq1}\frac{1}{|y|^{m+2\gamma
_{0}-2\alpha}}\int_{\tau}^{\tau+T}\int_{\mathbb{R}^{m}}\left(
u(t,x+y)+u(t,x-y)-2u(t,x)\right)  ^{p}dxdtdy\\[1ex]
&  \leq K_{4}\left\Vert u\right\Vert _{L^{p}((\tau,\tau+T)\times\mathbb{R}%
^{m})}^{p}.
\end{split}
\label{eq36}%
\end{equation}
Thus, Lebesgue's theorem, together with \eqref{eq35}-\eqref{eq36}, concludes
the proof of this lemma.
\end{proof}

\bigskip

In the same way, we can prove the following result.

\begin{lemma}
\label{ConvergLaplacianGamma2}For any $u\in C_{0}^{\infty}(\mathbb{R}^{m})$,
$\lim_{\gamma\rightarrow\gamma_{0}^{-}}(-\Delta)^{\gamma}u=\left(
-\Delta\right)  ^{\gamma_{0}}u$ strongly in $L^{p}(\mathbb{R}^{m})$ for any
$p\geq1$ and $\gamma_{0}\in\left(  0,1\right)  $.
\end{lemma}

Our objective now is to obtain some properties of $(-\Delta)^{\frac{\gamma}%
{2}} u$, $\gamma\in(0,1)$, regarding the fractional Laplacian operators. Let
\[
D_{\gamma}=D(\left(  -\Delta\right)  ^{\gamma})=\{u\in L^{2}(\mathbb{R}%
^{m}):\left(  -\Delta\right)  ^{\gamma}u\in L^{2}(\mathbb{R}^{m})\}.
\]

\begin{lemma}
\label{DinDs}For any $u\in H^{2}(\mathbb{R}^{m})$ and $\gamma\in\left(
0,1\right)  $, there exists a constant $C_{\gamma}>0$ such that%
\[
\left\Vert \left(  -\Delta\right)  ^{\gamma}u\right\Vert \leq C_{\gamma
}\left\Vert u\right\Vert _{H^{2}(\mathbb{R}^{m})}.
\]
In particular, the embedding $H^{2}(\mathbb{R}^{m})\subset D_{\gamma}$ is continuous.
\end{lemma}

\begin{proof}
For $u\in H^{2}(\mathbb{R}^{m})$, let $u_{n}\in C_{0}^{\infty}\left(
\mathbb{R}^{m}\right)  $ be a sequence converging to $u$ in $H^{2}%
(\mathbb{R}^{m})$. Observe the H\"older inequality implies that
\begin{equation}%
\begin{split}
&  ~\quad\frac{1}{4}C^{2}(m,\gamma)\int_{\mathbb{R}^{m}}\left\vert
\int_{\left\vert y\right\vert <1}\frac{u(x+y)+u(x-y)-2u(x)}{|y|^{m+2\gamma}%
}dy\right\vert ^{2}dx\\
&  \leq\frac{1}{4}C^{2}(m,\gamma)\int_{\left\vert y\right\vert <1}\frac
{1}{|y|^{m+2(\gamma-1)}}dy\int_{\mathbb{R}^{m}}\int_{\left\vert y\right\vert
<1}\frac{\left(  u(x+y)+u(x-y)-2u(x)\right)  ^{2}}{|y|^{m+2\gamma+2}}dydx.
\end{split}
\label{EmbEst1}%
\end{equation}
If the last integral is well defined and can be bounded by a suitable estimate
in terms of the norm $\left\Vert u\right\Vert _{H^{2}(\mathbb{R}^{m})}$, the
result of this lemma holds.

First, it is clear for a.a. $x\in\mathbb{R}^{m}$ and $y\in B_{1}$ that
\begin{equation}
\frac{\left(  u_{n}(x+y)+u_{n}(x-y)-2u_{n}(x)\right)  ^{2}}{|y|^{m+2\gamma+2}%
}\overset{n\rightarrow\infty}{\longrightarrow}\frac{\left(
u(x+y)+u(x-y)-2u(x)\right)  ^{2}}{|y|^{m+2\gamma+2}}. \label{eq38}%
\end{equation}
On the one hand, making use of a similar idea as in Lemma
\ref{ConvergLaplacianGamma}, we infer that there exist some positive constants
{$C_{1}$} and $C_{2,\gamma}$ such that
\[%
\begin{split}
&  ~\quad\int_{\mathbb{R}^{m}}\int_{\left\vert y\right\vert <1}\frac{\left(
u_{n}(x+y)+u_{n}(x-y)-2u_{n}(x)\right)  ^{2}}{|y|^{m+2\gamma+2}}dydx\\
&  \leq\int_{\mathbb{R}^{m}}\int_{\left\vert y\right\vert <1}\frac
{1}{|y|^{m+2\gamma+2}}\left(  \int_{0}^{1}\int_{-r}^{r}D^{2}u_{n}%
(x+sy)y\text{\textperiodcentered}ydsdr\right)  ^{2}dydx\\
&  \leq\int_{\mathbb{R}^{m}}\int_{\left\vert y\right\vert <1}\frac
{1}{|y|^{m+2(\gamma-1)}}\left(  \int_{-1}^{1}\left\vert D^{2}u_{n}%
(x+sy)\right\vert ds\right)  ^{2}dydx\\
&  \leq C_{1}\int_{\left\vert y\right\vert <1}\frac{1}{|y|^{m+2\left(
\gamma-1\right)  }}\int_{-1}^{1}\int_{\mathbb{R}^{m}}\left\vert D^{2}%
u_{n}(x+sy)\right\vert ^{2}dxdsdy\\
&  \leq C_{1}\int_{\left\vert y\right\vert <1}\frac{1}{|y|^{m+2\left(
\gamma-1\right)  }}\int_{-1}^{1}\int_{\mathbb{R}^{m}}\left\vert D^{2}%
u_{n}(x)\right\vert ^{2}dxdsdy\\[1ex]
&  \leq C_{2,\gamma}\int_{\mathbb{R}^{m}}|D^{2}u_{n}(x)|^{2}dx.
\end{split}
\]
Then, \eqref{eq38} and Fatou's lemma imply that there exist positive constants
$C_{3,\gamma}$ such that%
\begin{equation}%
\begin{split}
&  ~\quad\int_{\mathbb{R}^{m}}\int_{\left\vert y\right\vert <1}\frac{\left(
u(x+y)+u(x-y)-2u(x)\right)  ^{2}}{|y|^{m+2\gamma+2}}dydx\\
&  \leq\liminf_{n\rightarrow\infty}\int_{\mathbb{R}^{m}}\int_{\left\vert
y\right\vert <1}\frac{\left(  u_{n}(x+y)+u_{n}(x-y)-2u_{n}(x)\right)  ^{2}%
}{|y|^{m+2\gamma+2}}dydx\\
&  \leq C_{2,\gamma}\lim_{n\rightarrow\infty}\int_{\mathbb{R}^{m}}\left\vert
D^{2}u_{n}(x)\right\vert ^{2}dx\\
&  =C_{2,\gamma}\int_{\mathbb{R}^{m}}\left\vert D^{2}u(x)\right\vert
^{2}dx\leq{C_{3,\gamma}}\left\Vert u\right\Vert _{H^{2}(\mathbb{R}^{m})}^{2}.
\end{split}
\label{EmbEst2}%
\end{equation}
On the other hand, by the H\"{o}lder inequality, we know there exist some
positive constants $C_{4,\gamma}$ and $C_{5,\gamma}$ such that
\begin{equation}%
\begin{split}
&  \quad~\frac{1}{4}C^{2}(m,\gamma)\int_{\mathbb{R}^{m}}\left\vert
\int_{\left\vert y\right\vert \geq1}\frac{u(x+y)+u(x-y)-2u(x)}{|y|^{m+2\gamma
}}dy\right\vert ^{2}dx\\
&  \leq\frac{1}{4}C^{2}(m,\gamma)\int_{\mathbb{R}^{m}}\int_{\left\vert
y\right\vert \geq1}\frac{1}{|y|^{m+2\gamma}}dy\int_{\left\vert y\right\vert
\geq1}\frac{\left(  u(x+y)+u(x-y)-2u(x)\right)  ^{2}}{|y|^{m+2\gamma}}dydx\\
&  \leq C_{4,\gamma}\int_{\left\vert y\right\vert \geq1}\frac{1}%
{|y|^{m+2\gamma}}\int_{\mathbb{R}^{m}}\left(  u(x+y)+u(x-y)-2u(x)\right)
^{2}dxdy\\[1ex]
&  \leq C_{5,\gamma}\left\Vert u\right\Vert ^{2}.
\end{split}
\label{EmbEst3}%
\end{equation}
Collecting together (\ref{EmbEst1})-(\ref{EmbEst3}), we obtain that $\left(
-\Delta\right)  ^{\gamma}u\in L^{2}(\mathbb{R}^{m})$ and
\[
\left\Vert \left(  -\Delta\right)  ^{\gamma}u\right\Vert \leq C_{\gamma
}\left\Vert u\right\Vert _{H^{2}(\mathbb{R}^{m})},
\]
for some constant $C_{\gamma}:=\max\{C_{3,\gamma},C_{5,\gamma}\}$. The proof
of this lemma is complete.
\end{proof}

\begin{lemma}
\label{IntegrParts} If $u\in C_{0}^{\infty}(\mathbb{R}^{m})$, then%
\[
\left(  \left(  -\Delta\right)  ^{\frac{\gamma}{2}}u,\left(  -\Delta\right)
^{\frac{\gamma}{2}}u\right)  =\left(  \left(  -\Delta\right)  ^{\gamma
}u,u\right)  , \quad\text{ for any }\gamma\in\left(  0,1\right)  .
\]

\end{lemma}

\begin{proof}
By Proposition 3.6 in \cite{Nezza}, for any $u\in V_{\gamma}$, we have%
\[
\left(  \left(  -\Delta\right)  ^{\frac{\gamma}{2}}u,\left(  -\Delta\right)
^{\frac{\gamma}{2}}u\right)  =\left\Vert \left(  -\Delta\right)
^{\frac{\gamma}{2}}u\right\Vert ^{2}=\frac{1}{2}C(m,s)\left\Vert u\right\Vert
_{\dot{H}^{\gamma}(\mathbb{R}^{m})}^{2}.
\]
Meanwhile, it follows from \eqref{eq23} that%
\[%
\begin{split}
\frac{1}{2}C(m,\gamma)\Vert u\Vert_{\dot{H}^{\gamma}(\mathbb{R}^{m})}^{2}  &
=\frac{1}{2}C(m,\gamma)\int_{\mathbb{R}^{m}}\int_{\mathbb{R}^{m}}%
\frac{|u(x)-u(y)|^{2}}{|x-y|^{m+2\gamma}}dxdy\\[1ex]
&  =\left\langle A^{\gamma}(u),u\right\rangle _{(V_{\gamma}^{\ast},V_{\gamma
})}.
\end{split}
\]
Hence, for $u\in C_{0}^{\infty}(\mathbb{R}^{m})$, we have%
\[%
\begin{split}
&  \quad~\left(  \left(  -\Delta\right)  ^{\frac{\gamma}{2}}u,\left(
-\Delta\right)  ^{\frac{\gamma}{2}}u\right) \\
&  =\frac{1}{2}C(m,\gamma)\int_{\mathbb{R}^{m}}\int_{\mathbb{R}^{m}}%
\frac{|u(x)-u(y)|^{2}}{|x-y|^{m+2\gamma}}dxdy\\
&  =C(m,\gamma)\int_{\mathbb{R}^{m}}u(x)P.V.\int_{\mathbb{R}^{m}}%
\frac{u(x)-u(y)}{|x-y|^{m+2\gamma}}dydx\\[1ex]
&  =\left(  \left(  -\Delta\right)  ^{\gamma}u,u\right)  .
\end{split}
\]
Here $P.V.$ is a commonly used abbreviation for \textquotedblleft in the
principal value sense", for more details, see \cite{Nezza}.
\end{proof}

\begin{lemma}
\label{Norm1} If $u\in C_{0}^{\infty}(\mathbb{R}^{m})$, then $\left\Vert
\left(  -\Delta\right)  ^{\frac{1}{2}}u\right\Vert =\left\Vert \nabla
u\right\Vert .$
\end{lemma}

\begin{proof}
Taking a sequence $\alpha_{n}\rightarrow1^{-}$ as $n\rightarrow\infty$, making
use of lemmas \ref{ConvergLaplacianLp2} and \ref{IntegrParts}, we obtain
\[
\left\Vert \left(  -\Delta\right)  ^{\frac{\alpha_{n}}{2}}u\right\Vert
^{2}=\left(  \left(  -\Delta\right)  ^{\frac{\alpha_{n}}{2}}u,\left(
-\Delta\right)  ^{\frac{\alpha_{n}}{2}}u\right)  =\left(  \left(
-\Delta\right)  ^{\alpha_{n}}u,u\right)
\xrightarrow[n\rightarrow \infty]{\alpha_n\rightarrow 1^-}\left(  -\Delta
u,u\right)  =\left\Vert \nabla u\right\Vert ^{2}.
\]
Moreover, it follows from Lemma \ref{ConvergLaplacianGamma2} that,%
\[
\left\Vert \left(  -\Delta\right)  ^{\frac{\alpha_{n}}{2}}u\right\Vert
\xrightarrow[n\rightarrow \infty]{\alpha_n\rightarrow 1^-}\left\Vert \left(
-\Delta\right)  ^{\frac{1}{2}}u\right\Vert .
\]
We conclude the proof by the uniqueness of the limit.
\end{proof}

\begin{lemma}
\label{Norm2}If $u\in H^{2}(\mathbb{R}^{m})$, then $\left\Vert \left(
-\Delta\right)  ^{\frac{1}{2}}u\right\Vert =\left\Vert \nabla u\right\Vert .$
\end{lemma}

\begin{proof}
We take a sequence $\{u_{n}\}\subset C_{0}^{\infty}(\mathbb{R}^{m})$
converging to $u$ in $H^{2}(\mathbb{R}^{m})$. First, it is easy to see from
Lemma \ref{Norm1} that
\[
\left\Vert \left(  -\Delta\right)  ^{\frac{1}{2}}u_{n}\right\Vert =\left\Vert
\nabla u_{n}\right\Vert .
\]
Next, as $u_{n}\rightarrow u$ in $V$, we immediately derive that
\[
\left\Vert \nabla u_{n}\right\Vert \overset{n\rightarrow\infty
}{\longrightarrow}\left\Vert \nabla u\right\Vert .
\]
At last, we deduce from Lemma \ref{DinDs} that
\[
\left\Vert \left(  -\Delta\right)  ^{\frac{1}{2}}u_{n}\right\Vert
\overset{n\rightarrow\infty}{\longrightarrow}\left\Vert \left(  -\Delta
\right)  ^{\frac{1}{2}}u\right\Vert .
\]
Therefore, $\left\Vert \left(  -\Delta\right)  ^{\frac{1}{2}}u\right\Vert
=\left\Vert \nabla u\right\Vert $.
\end{proof}

\section{Convergence of solutions}

\label{s4} This section is devoted to proving the solutions of problem
(\ref{Eq0}) converge as $\gamma\rightarrow1^{-}$ to the solutions of the limit
problem with $\gamma=1$, that is, to the solution of the standard
reaction-diffusion equation \eqref{Eq1}. To that end, we begin by proving the
hemicontinuity of function $f$.

For each $t\in\mathbb{R}$, let $F(t,\cdot):H\rightarrow H$ be the Nemytskii
operator defined by $y=F(t,u)$, if $y\left(  x\right)
=F(t,u)(x):=f(t,x,u(x))$ for a.a. $x\in\mathbb{R}^{m}$. {In view of
(\ref{Growth}), we know that $y\in H$ if $u\in H$.} Hence, this operator is
well defined.

\begin{lemma}
\label{Hemicont} For each $t\in\mathbb{R}$, the function $F(t,\cdot)$ is
hemicontinuous, i.e., the real map $\lambda\longmapsto\left(  F(t, u+\lambda
v),w\right)  $ is continuous for any $u,v,w\in H.$
\end{lemma}

\begin{proof}
It is obvious that for a.a. $x\in\mathbb{R}^{m}$,
\[
f(t,x,u(x)+\lambda v(x))w(x)\rightarrow f(t,x,u(x)+\lambda_{0}%
v(x))w(x),~~\text{as }\lambda\rightarrow\lambda_{0}\in\mathbb{R}.
\]
By (\ref{Growth}), we know that
\[%
\begin{split}
\left\vert f(t,x,u(x)+\lambda v(x))w(x)\right\vert  &  \leq\left(  \psi
_{3}(t,x)+\psi_{2}(t,x)\left(  \left\vert u(x)\right\vert +\lambda\left\vert
v(x)\right\vert \right)  \right)  \left\vert w(x)\right\vert \\[0.8ex]
&  =r(t,x),~~\text{ for a.a. }x\in\mathbb{R}^{m},
\end{split}
\]
where $r(t,$\textperiodcentered$)\in L^{1}(\mathbb{R}^{m})$. Thus, the
Lebesgue theorem implies that%
\[%
\begin{split}
\left(  F(t,u+\lambda v),w\right)   &  =\int_{\mathbb{R}^{m}}%
f(t,x,u(x)+\lambda v(x))w(x)dx\\
&  \overset{\lambda\rightarrow\lambda_{0}}{\longrightarrow}\int_{\mathbb{R}%
^{m}}f(t,x,u(x)+\lambda_{0}v(x))w(x)dx=\left(  F(t,u+\lambda_{0}v),w\right)  .
\end{split}
\]
Therefore, $F(t,$\textperiodcentered$)$ is hemicontinuous.
\end{proof}

\begin{theorem}
\label{ConvergSolutions}Let $u_{\tau}^{n}\rightarrow u_{\tau}$ in $H$ {as
$n\rightarrow\infty$,} and let $u_{n}(\cdot)$ be the solution of problem
(\ref{Eq0}) for $\gamma=\gamma_{n}\in\left(  0,1\right)  $ with initial value
$u_{\tau}^{n}$. If $\gamma_{n}\rightarrow1^{-}$ {as $n\rightarrow\infty$},
then, $u_{n}\rightarrow u$ weak-star in $L^{\infty}(\tau,\tau+T;H)$ and weakly
in $L^{2}(\tau,\tau+T;H)$ as $n\rightarrow\infty$ for any $T>0$, where
$u\left(  \text{\textperiodcentered}\right)  $ is the unique solution of
problem (\ref{Eq0}) with $\gamma=1$. Moreover, for any sequence $t_{n}%
\in\lbrack\tau,\tau+T]$ converging to $t_{0}\in\lbrack\tau,\tau+T]$, we have%
\[
{u_{n}(t_{n})\rightarrow u\left(  t_{0}\right)  \text{ weakly in }H}.
\]

\end{theorem}

\begin{proof}
Assume first that $\sigma=0$ in (\ref{Lip}).

Let $\gamma_{n}\in(0,1)$ and $\gamma_{n}\rightarrow1^{-}$ as $n\rightarrow
\infty$. Denote by $u_{n}\left(  \text{\textperiodcentered}\right)  $ the
unique solution to problem (\ref{Eq0}) with initial condition $u_{\tau}^{n}$.
Multiplying equation (\ref{Eq0}) by $u_{n}$, making use of (\ref{Diss}) and
the Young inequality, we have%
\begin{equation}%
\begin{split}
\frac{1}{2}\frac{d}{dt}\left\Vert u_{n}\right\Vert ^{2}+\mu\left\Vert
u_{n}\right\Vert ^{2}+\frac{1}{2}C(m,\gamma_{n})\Vert u_{n}\Vert_{\dot
{H}^{\gamma_{n}}(\mathbb{R}^{m})}^{2}  &  \leq\int_{\mathbb{R}^{m}}\left(
\psi_{1}(t,x)+h(t,x)\right)  u_{n}(t,x)dx\\
&  \leq\int_{\mathbb{R}^{m}}\psi_{1}(t,x)dx+\frac{1}{2\mu}\left\Vert
h(t)\right\Vert ^{2}+\frac{\mu}{2}\left\Vert u_{n}\right\Vert ^{2}.
\end{split}
\label{Ineq1}%
\end{equation}
Consequently,
\[
\frac{d}{dt}\left\Vert u_{n}\right\Vert ^{2}+\mu\left\Vert u_{n}\right\Vert
^{2}+C(m,\gamma_{n})\Vert u_{n}\Vert_{\dot{H}^{\gamma_{n}}(\mathbb{R}^{m}%
)}^{2}\leq2\int_{\mathbb{R}^{m}}\psi_{1}(t,x)dx+\frac{1}{\mu}\left\Vert
h(t)\right\Vert ^{2}.
\]
Multiplying the above inequality by $e^{\mu s}$ and integrating it over
$\left(  \tau,t\right)  $ with $t>\tau$, we obtain%
\begin{equation}
\left\Vert u_{n}(t)\right\Vert ^{2}\leq\left\Vert u_{\tau}^{n}\right\Vert
^{2}e^{-\mu(t-\tau)}+2\int_{\tau}^{t}e^{-\mu(t-s)}\int_{\mathbb{R}^{m}}%
\psi_{1}(s,x)dxds+\frac{1}{\mu}\int_{\tau}^{t}e^{-\mu(t-s)}\left\Vert
h(s)\right\Vert ^{2}ds. \label{Ineq2}%
\end{equation}
Subsequently, for $T>0$, there exist two positive constants $M_{T}$ (which
depend on $T$) and $\overline{M}_{T,\gamma_{n}}$ (which depend on $T$ and
$\gamma_{n}$), such that%
\begin{equation}
\sup_{t\in\lbrack\tau,\tau+T]}\ \left\Vert u_{n}(t)\right\Vert \leq M_{T},
\label{Est1}%
\end{equation}
and%
\begin{equation}
C(m,\gamma_{n})\int_{\tau}^{\tau+T}\Vert u_{n}(s)\Vert_{\dot{H}^{\gamma_{n}%
}(\mathbb{R}^{m})}^{2}ds\leq\overline{M}_{T,\gamma_{n}}. \label{Est2}%
\end{equation}
Thus, $\{u_{n}\}$ is bounded in $L^{\infty}(\tau,\tau+T;H)$. It also follows
from \eqref{Est2} that $\{(-\Delta)^{\frac{\gamma_{n}}{2}}u_{n}\}$ is bounded
in $L^{2}(\tau,\tau+T;H)$, therefore $\{A^{\gamma_{n}}(u_{n})\}$ is bounded in
$L^{2}(\tau,\tau+T;V^{\ast})$. Moreover, by (\ref{Est1}) and (\ref{Growth}),
we infer that $\{F($\textperiodcentered$,u_{n})\}$ is bounded in $L^{2}%
(\tau,\tau+T;H)$. Hence, $\{\frac{du_{n}}{dt}\}$ is bounded in $L^{2}%
(\tau,\tau+T;V^{\ast})$. Then there exist functions $u,\chi$ and $\zeta$ such
that, up to a subsequence (relabeled the same), the following convergences
take place,
\begin{equation}
u_{n}\rightarrow u\text{ weak-star in }L^{\infty}(\tau,\tau+T;H),
\label{Converg1}%
\end{equation}
\begin{equation}
u_{n}\rightarrow u\text{ weakly in }L^{2}(\tau,\tau+T;H), \label{Converg3}%
\end{equation}%
\begin{equation}
F(\text{\textperiodcentered},u_{n})\rightarrow\chi\text{ weakly in }L^{2}%
(\tau,\tau+T;H), \label{Converg4}%
\end{equation}%
\begin{equation}
A^{\gamma_{n}}(u_{n})\rightarrow\zeta\text{ weakly in }L^{2}(\tau
,\tau+T;V^{\ast}), \label{Converg5}%
\end{equation}%
\begin{equation}
\frac{du_{n}}{dt}\rightarrow\frac{du}{dt}\text{ weakly in }L^{2}(\tau
,\tau+T;V^{\ast}). \label{Converg6}%
\end{equation}

Let us first check $\zeta=-\Delta u$. To this end, for any $\varphi\in
C_{0}^{\infty}(\left(  \tau,\tau+T\right)  \times\mathbb{R}^{m})$, by Lemma
\ref{ConvergLaplacianLp}, we have%
\[%
\begin{split}
&  \quad~<A^{\gamma_{n}}(u_{n}),\varphi>_{(\mathcal{D}^{\prime},\mathcal{D}%
)}=\int_{\tau}^{\tau+T}<A^{\gamma_{n}}(u_{n}),\varphi>_{(V^{\ast},V)}dt\\
&  =\int_{\tau}^{\tau+T}<u_{n},A^{\gamma_{n}}(\varphi)>_{(V,V^{\ast})}%
dt=\int_{\tau}^{\tau+T}\int_{\mathbb{R}^{m}}u_{n}(-\Delta)^{\gamma_{n}}\varphi
dxdt\\
&  \overset{\gamma_{n}\rightarrow1^{-}}{\longrightarrow}\int_{\tau}^{\tau
+T}\int_{\mathbb{R}^{m}}u(-\Delta)\varphi dxdt=\int_{\tau}^{\tau+T}%
<u,(-\Delta)\varphi>_{(V,V^{\ast})}dt=<(-\Delta)u,\varphi>_{(\mathcal{D}%
^{\prime},\mathcal{D})},
\end{split}
\]
where the above convergence follows from the facts that $u_{n}\rightarrow u$
weakly in $L^{2}(\tau,\tau+T;H)$ and $(-\Delta)^{\gamma_{n}}\varphi
\rightarrow(-\Delta)\varphi$ in $L^{2}(\tau,\tau+T;H)$. Therefore,
$A^{\gamma_{n}}(u_{n})\rightarrow-\Delta u$ as $\gamma^{n}\rightarrow1^{-}$ in
the sense of distributions, which implies that $\zeta=-\Delta u.$ Therefore,
for any $\eta\in L^{2}(\tau,\tau+T;V)$, the above convergences (\ref{Converg1}%
)-(\ref{Converg6}) imply that%
\[%
\begin{split}
&  \int_{\tau}^{\tau+T}<\frac{du_{n}}{dt},\eta>_{(V^{\ast},V)}dt+\int_{\tau
}^{\tau+T}<A^{\gamma_{n}}(u_{n}),\eta>_{(V^{\ast},V)}dt-\int_{\tau}^{\tau
+T}(F(t,u_{n})+h(t),\eta)dt\\[1.2ex]
&  \xrightarrow[n\rightarrow \infty]{\gamma_n\rightarrow 1^-}\int_{\tau}%
^{\tau+T}<\frac{du}{dt},\eta>_{(V^{\ast},V)}dt+\int_{\tau}^{\tau+T}<-\Delta
u,\eta>_{(V^{\ast},V)}dt-\int_{\tau}^{\tau+T}(\chi+h,\eta)dt=0.
\end{split}
\]
This equality is equivalent to (\ref{EqSolutionGamma}) (for more details on
this fact see e.g. \cite[p.43]{KMVY}).

In addition, we know that
\[
u\in L^{\infty}(\tau,\tau+T;H),
\]%
\[
\frac{du}{dt}\in L^{2}(\tau,\tau+T;V^{\ast}),
\]%
\[
-\Delta u\in L^{2}(\tau,\tau+T;V^{\ast}).
\]
Notice that, since $-\Delta u(t)\in V^{\ast}$ for a.a. $t\in(\tau,\tau+T)$, we
obtain that $u\left(  t\right)  \in V$. Therefore,%
\[
\left\Vert \nabla u(t)\right\Vert ^{2}=<-\Delta u(t),u(t)>_{(V^{\ast},V)}%
\leq\left\Vert \Delta u(t)\right\Vert _{V^{\ast}}\left\Vert \nabla
u(t)\right\Vert ,
\]
which implies,
\[
u\in L^{2}(\tau,\tau+T;V).
\]
Thus, $u$ is the {unique} weak solution of the problem%
\begin{equation}%
\begin{cases}
\dfrac{\partial u}{\partial t}-\Delta u=g(t),\quad(x,t)\in\mathbb{R}^{m}%
\times(\tau,\infty),\\
u(x,\tau)=u_{\tau}(x),\quad x\in\mathbb{R}^{m},
\end{cases}
\label{ProblemLimit}%
\end{equation}
where $g(t)=\chi(t)+h(t)\in L^{2}(\tau,\tau+T;H).$

Further, we need to prove that for any sequence $t_{n}\in\lbrack\tau,\tau+T]$
converging to $t_{0}\in\lbrack\tau,\tau+T]$, we have%
\[
{u_{n}(t_{n})\rightarrow u\left(  t_{0}\right)  \text{ weakly in }{H}}.
\]
It is clear that there is a subsequence $\{u_{n_{k}}(t_{n_{k}})\}$ and $y\in
H$ such that $u_{n_{k}}(t_{n_{k}})\rightarrow y$ weakly in $H$. If we prove
that $y=u(t_{0})$, then the result follows by a standard contradiction argument.

For any $R>0$, let us consider the spaces $H_{0}^{1}(B_{R})$, $L^{2}(B_{R})$,
$(H_{0}^{1}(B_{R}))^{\ast}$. Then we have the following chains of continuous
embeddings,%
\[
H_{0}^{1}(B_{R})\subset V\subset H\ \subset V^{\ast}\subset(H_{0}^{1}%
(B_{R}))^{\ast},
\]%
\[
H_{0}^{1}(B_{R})\subset L^{2}(B_{R})\subset(H_{0}^{1}(B_{R}))^{\ast},
\]
where we consider that $z\in H_{0}^{1}(B_{R})$ is an element of $V$ by setting
$z\left(  x\right)  =0$ for $x\not \in B_{R}$. Let $L_{R}z$ be the projection
of $z\in H$ onto $L^{2}(B_{R})$. Since the embedding $L^{2}(B_{R}%
)\subset(H_{0}^{1}(B_{R}))^{\ast}$ is compact, in view of estimate
(\ref{Ineq2}), we know that the sequence $\{L_{R}u_{n}(t)\}$ is relatively
compact in $(H_{0}^{1}(B_{R}))^{\ast}$ for any $t\in(\tau,\tau+T)$.
Furthermore, by means of the fact that $\{\frac{du_{n}(t)}{dt}\}$ is bounded
in $L^{2}(\tau,\tau+T;V^{\ast})\subset L^{2}(\tau,\tau+T;(H_{0}^{1}%
(B_{R}))^{\ast})$, we infer there exists a constant $C$, such that%
\[
\left\Vert u_{n}(t)-u_{n}(s)\right\Vert _{\left(  H_{0}^{1}(B_{R})\right)
^{\ast}}\leq\int_{s}^{t}\left\Vert \dfrac{du_{n}}{dt}(r)\right\Vert _{\left(
H_{0}^{1}(B_{R})\right)  ^{\ast}}dr\leq C\left(  t-s\right)  ^{\frac{1}{2}}.
\]
Hence, the sequence $\{L_{R}u_{n}(\cdot)\}$ is equicontinuous in $(H_{0}%
^{1}(B_{R}))^{\ast}$. The Ascoli-Arzel\`{a} theorem shows that, up to a
subsequence (relabeled the same), $L_{R}u_{n}\rightarrow L_{R}u$ in
$C([\tau,\tau+T],(H_{0}^{1}(B_{R}))^{\ast})$. By the compact embedding
$L^{2}(B_{R})\subset(H_{0}^{1}(B_{R}))^{\ast}$ and a diagonal argument, we
deduce that $L_{R}u_{n_{k}}(t_{n_{k}})\rightarrow L_{R}u(t_{0}),\ L_{R}%
u_{n_{k}}(t_{n_{k}})\rightarrow L_{R}y$ in $(H_{0}^{1}(B_{R}))^{\ast}$ for any
$R>0$. Therefore, $L_{R}u(t_{0})=L_{R}y$ in $(H_{0}^{1}(B_{R}))^{\ast}$ for
any $R>0$. We obtain then immediately that $u(t_{0})$ and $y$ are the same
distribution, hence, $y=u(t_{0})$.

At last, let us prove that $\chi=F($\textperiodcentered$,u($%
\textperiodcentered$))$. If $u_{\tau}\in V$, it is known (see Lemma
\ref{RegularityParabolic} below) that the unique solution $u\left(
\text{\textperiodcentered}\right)  $ to problem (\ref{ProblemLimit}) belongs
to {$L^{2}(\tau,\tau+T;H^{2}(\mathbb{R}^{m}))\cap L^{\infty}(\tau,\tau+T;V)$}.
As the unique weak solution $u\left(  \text{\textperiodcentered}\right)  $ of
(\ref{ProblemLimit}) satisfies that $u\left(  t\right)  \in V$ for a.a.
$t\in(\tau,\tau+T)$, it is easy to see that%
\begin{equation}
u\in L^{2}(\tau+\delta,\tau+T;H^{2}(\mathbb{R}^{m})),\quad\forall\delta
\in(0,T). \label{Regularity}%
\end{equation}
Now multiplying (\ref{Eq0}) by $u_{n}$, we obtain%
\[
\frac{1}{2}\frac{d}{dt}\left\Vert u_{n}\right\Vert ^{2}+\frac{1}{2}%
C(m,\gamma_{n})\Vert u_{n}\Vert_{\dot{H}^{\gamma_{n}}(\mathbb{R}^{m})}%
^{2}=\left(  F(t,u_{n}(t))+h(t,x),u_{n}(t)\right)  .
\]
Then, integrating the above equality over $(\tau,\tau+T)$, we find
\begin{equation}%
\begin{split}
\int_{\tau}^{\tau+T}\left(  F(t,u_{n}(t)),u_{n}(t)\right)  dt  &  =\frac{1}%
{2}\left\Vert u_{n}(\tau+T)\right\Vert ^{2}-\frac{1}{2}\left\Vert u_{n}%
(\tau)\right\Vert ^{2}\\
&  ~~+\frac{1}{2}C(m,\gamma_{n})\int_{\tau}^{\tau+T}\Vert u_{n}(t)\Vert
_{\dot{H}^{\gamma_{n}}(\mathbb{R}^{m})}^{2}dt\\
&  ~~-\int_{\tau}^{\tau+T}\left(  h(t,x),u_{n}(t)\right)  dt.
\end{split}
\label{eq412}%
\end{equation}
For every $v\in L^{2}(\tau,\tau+T;H)$, define the sequence%
\[
X_{n}:=\int_{\tau}^{\tau+T}\left(  F(t,u_{n}(t))-F(t,v(t)),u_{n}%
(t)-v(t)\right)  dt.
\]
Therefore, it follows from \eqref{eq412} and the above equality that%
\[%
\begin{split}
X_{n}  &  =\frac{1}{2}\left\Vert u_{n}(\tau+T)\right\Vert ^{2}-\frac{1}%
{2}\left\Vert u_{\tau}^{n}\right\Vert ^{2}+\frac{1}{2}C(m,\gamma_{n}%
)\int_{\tau}^{\tau+T}\Vert u_{n}(t)\Vert_{\dot{H}^{\gamma_{n}}(\mathbb{R}%
^{m})}^{2}dt\\[1ex]
&  ~~-\int_{\tau}^{\tau+T}\left(  h(t),u_{n}(t)\right)  dt-\int_{\tau}%
^{\tau+T}\left(  F(t,u_{n}(t)),v(t)\right)  dt-\int_{\tau}^{\tau+T}\left(
F(t,v(t)),u_{n}(t)-v(t)\right)  dt.
\end{split}
\]
On the one hand, by (\ref{Lip}), (\ref{Converg3}), (\ref{Converg4}), $u_{\tau
}^{n}\rightarrow u_{\tau}$ in $H$ and $u_{n}(\tau+T)\rightarrow u(\tau+T)$
weakly in $H$ as $n\rightarrow\infty$,  we obtain%
\begin{equation}%
\begin{split}
0  &  \geq\liminf_{n\rightarrow\infty}\ X_{n}\geq\frac{1}{2}\left\Vert
u(\tau+T)\right\Vert ^{2}-\frac{1}{2}\left\Vert u_{\tau}\right\Vert ^{2}%
+\frac{1}{2}\liminf_{n\rightarrow\infty}C(m,\gamma_{n})\int_{\tau}^{\tau
+T}\Vert u_{n}(t)\Vert_{\dot{H}^{\gamma_{n}}(\mathbb{R}^{m})}^{2}dt\\
&  ~~-\int_{\tau}^{\tau+T}\left(  h(t),u(t)\right)  dt-\int_{\tau}^{\tau
+T}\left(  \chi(t),v(t)\right)  dt-\int_{\tau}^{\tau+T}\left(
F(t,v(t)),u(t)-v(t)\right)  dt.
\end{split}
\label{eq413}%
\end{equation}
On the other hand, multiplying (\ref{ProblemLimit}) by $u$ and integrating it
over $(\tau,\tau+T)$, we have%
\begin{equation}
\frac{1}{2}\left\Vert u(\tau+T)\right\Vert ^{2}-\frac{1}{2}\left\Vert u_{\tau
}\right\Vert ^{2}-\int_{\tau}^{\tau+T}\left(  h(t),u(t)\right)  dt=\int_{\tau
}^{\tau+T}\left(  \chi(t),u(t)\right)  dt-\frac{1}{2}\int_{\tau}^{\tau
+T}\left\Vert \nabla u\right\Vert ^{2}dt. \label{eq414}%
\end{equation}
Thus, combing \eqref{eq413} with \eqref{eq414}, we deduce
\begin{equation}%
\begin{split}
0  &  \geq\int_{\tau}^{\tau+T}\left(  \chi(t)-F(t,v(t)),u(t)-v(t)\right)
dt\\[1ex]
&  ~~+\frac{1}{2}\liminf_{n\rightarrow\infty}C(m,\gamma_{n})\int_{\tau}%
^{\tau+T}\Vert u_{n}(t)\Vert_{\dot{H}^{\gamma_{n}}(\mathbb{R}^{m})}%
^{2}dt-\frac{1}{2}\int_{\tau}^{\tau+T}\left\Vert \nabla u\right\Vert ^{2}dt.
\end{split}
\label{eq415}%
\end{equation}
Observe that, for any $\varepsilon>0$, there exists $\delta(\varepsilon)>0$
such that
\[
{\frac{1}{2}\int_{\tau}^{\tau+\delta}\left\Vert \nabla u\right\Vert ^{2}%
dt\leq\varepsilon.}%
\]
Thus, \eqref{eq415} and the above inequality imply that
\begin{equation}%
\begin{split}
\varepsilon &  \geq\int_{\tau}^{\tau+T}\left(  \chi
(t)-F(t,v(t)),u(t)-v(t)\right)  dt\\[0.8ex]
&  ~~+\frac{1}{2}\liminf_{n\rightarrow\infty}C(m,\gamma_{n})\int_{\tau+\delta
}^{\tau+T}\Vert u_{n}(t)\Vert_{\dot{H}^{\gamma_{n}}(\mathbb{R}^{m})}%
^{2}dt-\frac{1}{2}\int_{\tau+\delta}^{\tau+T}\left\Vert \nabla u\right\Vert
^{2}dt.
\end{split}
\label{eq416}%
\end{equation}
Next, we want to show that%
\begin{equation}
\int_{\tau+\delta}^{\tau+T}\left\Vert \nabla u\right\Vert ^{2}dt\leq
\liminf_{n\rightarrow\infty}C(m,\gamma_{n})\int_{\tau+\delta}^{\tau+T}\Vert
u_{n}(t)\Vert_{\dot{H}^{\gamma_{n}}(\mathbb{R}^{m})}^{2}dt. \label{eq422}%
\end{equation}
It follows from (\ref{Est2}) and (\ref{EquivNorm}) that, up to a subsequence,
as $\gamma_{n}\rightarrow1^{-}$ and $n\rightarrow\infty$,
\begin{equation}
\left(  -\Delta\right)  ^{\frac{\gamma_{n}}{2}}u_{n}\rightarrow y,~~\text{
weakly in }L^{2}(\tau,\tau+T;H). \label{eq417}%
\end{equation}
If we are able to prove that $\left\Vert y\left(  t\right)  \right\Vert
=\left\Vert \nabla u\left(  t\right)  \right\Vert $ for a.a. $t\in(\tau
+\delta,\tau+T)$, then (\ref{eq416}) follows. To this end, for any
$\vartheta\in C_{0}^{\infty}(\left(  \tau,\tau+T\right)  \times\mathbb{R}%
^{m})$, we have
\begin{equation}
\int_{\tau+\delta}^{\tau+T}\left(  \left(  -\Delta\right)  ^{\frac{\gamma_{n}%
}{2}}u_{n}(t),\vartheta(t)\right)  dt=\int_{\tau+\delta}^{\tau+T}\left(
u_{n}(t),\left(  -\Delta\right)  ^{\frac{\gamma_{n}}{2}}\vartheta(t)\right)
dt. \label{eq418}%
\end{equation}
On the one hand, by \eqref{eq417}, we derive
\begin{equation}
\int_{\tau+\delta}^{\tau+T}\left(  \left(  -\Delta\right)  ^{\frac{\gamma_{n}%
}{2}}u_{n}(t),\vartheta(t)\right)  dt\rightarrow\int_{\tau+\delta}^{\tau
+T}(y(t),\vartheta(t))dt,\quad\text{as }\gamma_{n}\rightarrow1^{-}%
(n\rightarrow\infty). \label{eq419}%
\end{equation}
On the other hand, we deduce from Lemma \ref{ConvergLaplacianGamma} that
\begin{equation}
\int_{\tau+\delta}^{\tau+T}\left(  u_{n}(t),(-\Delta)^{\frac{\gamma_{n}}{2}%
}\vartheta(t)\right)  dt\rightarrow\int_{\tau+\delta}^{\tau+T}\left(
u(t),(-\Delta)^{\frac{1}{2}}\vartheta(t)\right)  dt,\quad\text{as }\gamma
_{n}\rightarrow1^{-}(n\rightarrow\infty). \label{eq420}%
\end{equation}
Moreover, by (\ref{Regularity}), we know that $u\left(  t\right)  \in
H^{2}(\mathbb{R}^{m})$ for a.a. $t\in(\tau+\delta,\tau+T)$, which, together
with Lemma \ref{DinDs}, imply that $\left(  -\Delta\right)  ^{\frac{1}{2}%
}u\left(  t\right)  \in H$. Therefore,
\begin{equation}
\int_{\tau+\delta}^{\tau+T}\left(  u(t),\left(  -\Delta\right)  ^{\frac{1}{2}%
}\vartheta(t)\right)  dt=\int_{\tau+\delta}^{\tau+T}\left(  \left(
-\Delta\right)  ^{\frac{1}{2}}u(t),\vartheta(t)\right)  dt. \label{eq421}%
\end{equation}
Thus, \eqref{eq418}-\eqref{eq421} yield $y=\left(  -\Delta\right)  ^{\frac
{1}{2}}u$ in $L^{2}(\tau+\delta,\tau+T;H)$. Then Lemma \ref{Norm2} shows%
\[
\left\Vert \nabla u(t)\right\Vert =\left\Vert \left(  -\Delta\right)
^{\frac{1}{2}}u(t)\right\Vert =\left\Vert y(t)\right\Vert ,~~\text{ for a.a.
}t\in(\tau+\delta,\tau+T),
\]
as desired.

Thus, from \eqref{eq416} and \eqref{eq422}, we infer that
\[
\int_{\tau}^{\tau+T}\left(  \chi(t)-F(t,v(t)),u(t)-v(t)\right)  dt\leq
\varepsilon,\quad\text{ for any }\varepsilon>0,
\]
which implies,%
\[
\int_{\tau}^{\tau+T}\left(  \chi(t)-F(t,v(t)),u(t)-v(t)\right)  dt\leq0.
\]
Letting $v\left(  t\right)  =u(t)-\theta z\left(  t\right)  $, where
$\theta>0$ and $z\in L^{2}(\tau,\tau+T;H)$, we have
\[
\int_{\tau}^{\tau+T}\left(  \chi(t)-F(t,u(t)-\theta z\left(  t\right)
),z(t)\right)  dt\leq0.
\]
Since $F(t,$\textperiodcentered$)$ is hemicontinuous (see Lemma \ref{Hemicont}%
), passing to the limit as $\theta\rightarrow0$ and using the Lebesgue
theorem, we derive%
\[
\int_{\tau}^{\tau+T}\left(  \chi(t)-F(t,u(t)),z(t)\right)  dt\leq0.
\]
As $z$ is arbitrary, we deduce that $\chi=F($\textperiodcentered
$,u($\textperiodcentered$)).$

Therefore, $u\left(  \text{\textperiodcentered}\right)  $ is the weak solution
to problem (\ref{Eq0}) with $\gamma=1$. Since every converging subsequence has
the same limit, the whole sequence satisfies the above convergences.

If $\sigma>0$ in (\ref{Lip}), then by Lemma \ref{EquivSol} we consider the
sequence of solutions $v_{n}(t)=u_{n}(t)e^{-\sigma t}$ of problem (\ref{Eq0C})
which, by the above, converge to the unique solution $v\left(
\text{\textperiodcentered}\right)  $ to problem (\ref{Eq1C}) with $v\left(
\tau\right)  =u_{\tau}e^{-\sigma\tau}$ in the sense given in the statement.
Hence, by Lemma \ref{EquivSol2} it readily follows that the convergences are
true for the sequence $\{u_{n}($\textperiodcentered$)\}$ as well.
\end{proof}

\bigskip

By a standard contradiction argument, we deduce the following result.

\begin{corollary}
Let $u_{\tau}^{\gamma}\rightarrow u_{\tau}$ as $\gamma\rightarrow1^{-}$. Let
$u_{\gamma}($\textperiodcentered$)$ be the solution to problem (\ref{Eq0})
with initial value $u_{\tau}^{\gamma}$. Then $u_{\gamma}\rightarrow u$ as
$\gamma\rightarrow1^{-}$ in the sense given in Theorem \ref{ConvergSolutions},
where $u\left(  \text{\textperiodcentered}\right)  $ is the unique solution to
problem (\ref{Eq0}) with initial value $u_{\tau}$ and $\gamma=1$.
\end{corollary}

\begin{lemma}
\label{RegularityParabolic}If $u_{\tau}\in V$, then the unique solution
$u\left(  \text{\textperiodcentered}\right)  $ to problem (\ref{ProblemLimit})
belongs to {$L^{2}(\tau,\tau+T;H^{2}(\mathbb{R}^{m}))\cap L^{\infty}(\tau
,\tau+T;V).$}
\end{lemma}

\begin{proof}
Making use of the standard procedures (see, e.g., the proof of \cite[Theorem
5]{MorVal}), we can approximate the solution $u\left(
\text{\textperiodcentered}\right)  $ by solutions of the problem in a bounded
domain with Dirichlet boundary conditions. Precisely, let us consider the
following parabolic equation, for any $R>0$,%
\begin{equation}%
\begin{cases}
\dfrac{\partial u}{\partial t}-\Delta u=g(t),\quad(x,t)\in B_{R}\times
(\tau,\infty),\\
u=0,\quad(x,t)\in\partial B_{R}\times(\tau,\infty),\\
u(x,\tau)=u_{\tau}^{R}(x),\quad x\in B_{R},
\end{cases}
\label{ProblemBR}%
\end{equation}
where $u_{\tau}^{R}(x)=\psi_{R}(\left\vert x\right\vert )u_{\tau}(x)$, with
$\psi_{R}$ being a smooth function such that\text{ }%
\begin{equation*}
\psi_{R}(s)=
\begin{cases}
\begin{split}
1,~~ \qquad \qquad ~~&\text{ if }~0\leq s\leq R-1,\\
0\leq\xi(s)\leq1,~~&\text{ if }~R-1< s< R,\\
0,~~\qquad\qquad~~ &\text{ if }~s\geq R,
\end{split}
\end{cases}
\end{equation*}
and $\left\vert \psi_R^{\prime}(s)\right\vert ,\ \left\vert \psi_R^{\prime\prime
}(s)\right\vert \leq L$, for all $s\geq0,\ R>0$, and some positive constant
$L$. If $u_{\tau}\in H$, then problem (\ref{ProblemBR}) has a unique weak
solution $u_{R}\left(  \text{\textperiodcentered}\right)  .$ We extend the
function $u_{R}$ to the whole space by defining
\[
\overline{u}_{R}(x)=%
\begin{cases}%
\begin{split}
\psi_{R}(\left\vert x\right\vert )u_{R}(x),\quad &  \text{ if }x\in B_{R},\\
0,\qquad\qquad\qquad &  \text{ otherwise.}%
\end{split}
\end{cases}
\]
Then, using the same arguments as in \cite{MorVal}, we obtain
\begin{align*}
\overline{u}_{R}  &  \rightarrow u~\text{ weak-star in }~L^{\infty}(\tau
,\tau+T;H),\\
\overline{u}_{R}  &  \rightarrow u~\text{ weakly in }~L^{\infty}(\tau
,\tau+T;V),
\end{align*}
as $R\rightarrow+\infty$.

If $u_{\tau}\in V$, then $u_{R}\left(  \text{\textperiodcentered}\right)  \in
L^{\infty}(\tau,\tau+T;H_{0}^{1}(B_{R}))\cap L^{2}(\tau,\tau+T;H^{2}%
(B_{R})\cap H_{0}^{1}(B_{R}))$ (see, for example, \cite[p.70]{Temam}).
Multiplying the equation in (\ref{ProblemBR}) by $-\Delta u_{R}$, we obtain%
\[
\frac{d}{dt}\left\Vert u_{R}\right\Vert _{H_{0}^{1}(B_{R})}^{2}+\left\Vert
\Delta u_{R}\right\Vert _{L^{2}(B_{R})}^{2}\leq\left\Vert g(t)\right\Vert
_{L^{2}(B_{R})}^{2}.
\]
Therefore, there exists a constant $K>0$ such that,
\[
\left\Vert u_{R}(t)\right\Vert _{H_{0}^{1}(B_{R})}^{2}+\int_{\tau}^{\tau
+T}\left\Vert \Delta u_{R}\right\Vert _{L^{2}(B_{R})}^{2}ds\leq\left\Vert
u_{\tau}^{R}\right\Vert _{H_{0}^{1}(B_{R})}^{2}+\int_{\tau}^{\tau+T}\left\Vert
g(s)\right\Vert _{L^{2}(B_{R})}^{2}ds\leq K.
\]
It follows that $\overline{u}_{R}$ is bounded in {$L^{\infty}(\tau
,\tau+T;V)\cap L^{2}(\tau,\tau+T;H^{2}(\mathbb{R}^{m}))$. Hence,}%
\begin{align*}
\overline{u}_{R}  &  \rightarrow u~\text{ weak-star in }~L^{\infty}(\tau
,\tau+T;V),\\
\overline{u}_{R}  &  \rightarrow u~\text{ weakly in }~L^{2}(\tau,\tau
+T;H^{2}(\mathbb{R}^{m})),
\end{align*}
as $R\rightarrow+\infty$, and the result follows.
\end{proof}

\section{Existence of global attractors}

\label{s5}

We will focus on the existence of global attractors to problem \eqref{Eq0} in
this section for a more general nonlinear term $f$. To show the main ideas,
assume that the functions $h\in H$, $f:\mathcal{\mathbb{R}}^{m}\times
\mathcal{\mathbb{R}}\rightarrow\mathbb{R}$ do not depend on time. Namely, we
will study the autonomous case of \eqref{Eq0} and take $\tau=0$. In addition,
suppose that $f$ is a continuous function such that, for all $x\in
\mathbb{R}^{m}$ and $u\in\mathbb{R}$, $f(x,\cdot)\in C^{1}(\mathbb{R})$ and%
\begin{equation}
\frac{\partial f}{\partial u}(x,u)\leq\sigma,\label{LipAut}%
\end{equation}%
\begin{equation}
f(x,u)u\leq-\beta|u|^{p}+\psi_{1}(x),\label{DissAut}%
\end{equation}%
\begin{equation}
|f(x,u)|\leq\psi_{2}(x)|u|^{p-1}+\psi_{3}(x),\label{GrowthAut}%
\end{equation}
where $\beta>0$, $\sigma\geq0$, $p\geq2$ are constants, $\psi_{1}\in
L^{1}(\mathbb{R}^{m}),\ \psi_{2}\in L^{\infty}(\mathbb{R}^{m}),\ \psi_{3}\in
L^{q}(\mathbb{R}^{m})$ and $q$ is the conjugate of $p$, that is, $1/p+1/q=1.$

Define the Nemytskii operator $F:L^{p}(\mathbb{R}^{m})\rightarrow
L^{q}(\mathbb{R}^{m})$ by $y=F(u)$, if $y\left(  x\right)  =f(x,u(x))$ for
a.a. $x\in\mathbb{R}^{m}$. Now, we consider the problem%
\begin{equation}%
\begin{cases}
\dfrac{\partial u}{\partial t}+(-\Delta)^{\gamma}u+\mu u=f(x,u)+h(x),\quad
(x,t)\in\mathbb{R}^{m}\mathcal{\times}\left(  0,\infty\right)  ,\\
u(x,0)=u_{0}(x),\quad x\in\mathbb{R}^{m},
\end{cases}
\label{EqAut}%
\end{equation}
where $\mu>0.$

\begin{definition}
Let $u_{0}\in H$ and $\gamma\in(0,1)$. A continuous function $u:[0,\infty
)\rightarrow H$ is said to be a weak solution of problem (\ref{EqAut}) if
$u(0)=u_{0}$, $u\in L_{loc}^{2}(0,\infty;V_{\gamma})\cap L_{loc}^{p}%
(0,\infty;L^{p}(\mathbb{R}^{m}))$, $\dfrac{du}{dt}\in L_{loc}^{2}%
(0,\infty;V_{\gamma}^{\ast})+L_{loc}^{q}(0,\infty;L^{q}(\mathbb{R}^{m}))$, and
$u$ satisfies, for every $\xi\in V_{\gamma}\cap L^{p}(\mathbb{R}^{m})$,
\begin{align}
\frac{d}{dt}\left(  u,\xi\right)   &  +\frac{1}{2}C(m,\gamma)\int%
_{\mathbb{R}^{m}}\int_{\mathbb{R}^{m}}\frac{(u(x)-u(y))(\xi(x)-\xi
(y))}{|x-y|^{m+2\gamma}}dxdy+\mu(u,\xi)\label{EqualitySol}\\
&  =\int_{\mathbb{R}^{m}}\left(  f(x,u(t,x))+h(x)\right)  \xi\left(  x\right)
dx,\nonumber
\end{align}
in the sense of scalar distributions in $\left(  0,\infty\right)  .$
\end{definition}

It is known (see, for example, \cite[Theorem 3.2]{LuQiWangZhang}) that problem
(\ref{EqAut}) possesses a unique weak solution for any $u_{0}\in H$ and
$\gamma\in\left(  0,1\right)  $, which is continuous with respect to the
initial datum $u_{0}$. For any $\gamma\in(0,1)$, let us define the family of
operators $T_{\gamma}:\mathbb{R}^{+}\times H\rightarrow H$, given by%
\[
T_{\gamma}(t,u_{0})=u_{\gamma}(t),
\]
where $u_{\gamma}(\cdot) $ is the unique weak solution to problem
(\ref{EqAut}) with initial datum $u\left(  0\right)  =u_{0}$. It is not
difficult to see that the map $\left(  t,u_{0}\right)  \mapsto T_{\gamma
}(t,u_{0})$ is continuous.

Before showing the main results, let us first recall some definitions which
are essential for our proofs.

\begin{definition}
The compact set $\mathcal{A}_{\gamma}$ is said to be a global attractor for
$T_{\gamma}$, if the following conditions are fulfilled:

\begin{enumerate}
\item $T_{\gamma}(t,\mathcal{A}_{\gamma})=\mathcal{A}_{\gamma}$, for all
$t\geq0$ (invariance);

\item $\mathcal{A}_{\gamma}$ attracts any bounded set $B$ of $H$, that is:
\[
dist(T_{\gamma}(t,B),\mathcal{A}_{\gamma})\rightarrow0\text{ as }%
t\rightarrow+\infty,
\]
where $dist(C,D)=\sup_{x\in C}\inf_{y\in D}\ \left\Vert x-y\right\Vert $ is
the Hausdorff semidistance in $H.$
\end{enumerate}
\end{definition}

\begin{definition}
The map $\phi:\mathbb{R}\rightarrow H$ is called a complete trajectory of
$T_{\gamma}$, if $\phi(t)=T_{\gamma}(t-s,\phi(s))$ for any $t\geq s.$
\end{definition}

\begin{lemma}
\label{Estimates}Any solution $u\left(  \text{\textperiodcentered}\right)  $
of (\ref{EqAut}) satisfies%
\begin{equation}
\left\Vert u(t)\right\Vert ^{2}+C(m,\gamma)\int_{0}^{t}e^{-\mu(t-s)}\Vert
u\Vert_{\dot{H}^{\gamma}(\mathbb{R}^{m})}^{2}ds\leq\left\Vert u_{0}\right\Vert
^{2}e^{-\mu t}+\frac{2}{\mu}\int_{\mathbb{R}^{m}}\psi_{1}(x)dx+\frac{\Vert
h\Vert^{2}}{\mu^{2}}. \label{EstsolGamma}%
\end{equation}

\end{lemma}

\begin{proof}
Multiplying equation (\ref{EqAut}) by $u$, making using of (\ref{DissAut}) and
the Young inequality, we have%
\begin{equation}
\frac{d}{dt}\left\Vert u\right\Vert ^{2}+\mu\left\Vert u\right\Vert
^{2}+C(m,\gamma)\Vert u\Vert_{\dot{H}^{\gamma}(\mathbb{R}^{m})}^{2}%
+2\beta\left\Vert u\right\Vert _{p}^{p}\leq2\int_{\mathbb{R}^{m}}\psi
_{1}(x)dx+\frac{1}{\mu}\left\Vert h\right\Vert ^{2}. \label{EstSol2}%
\end{equation}
Multiplying the above inequality by $e^{\mu s}$ and integrating over $\left(
0,t\right)  $ with $t>0$, we obtain (\ref{EstsolGamma}).
\end{proof}

\begin{lemma}
\label{Absorbing}For any $R_{1}>0$, there exists $T(R_{1})>0$ such that,%
\[
\left\Vert T_{\gamma}(t,u_{0})\right\Vert \leq R_{0}:=\sqrt{1+\frac{2}{\mu
}\int_{\mathbb{R}^{m}}\psi_{1}(x)dx+\frac{\Vert h\Vert^{2}}{\mu^{2}}}%
,\qquad\forall t\geq T,
\]
for any $u_{0}$ satisfying $\left\Vert u_{0}\right\Vert \leq R_{1}.$
\end{lemma}

\begin{proof}
This lemma follows easily from (\ref{EstsolGamma}), so we omit the details of
the proof here.
\end{proof}

\begin{remark}
The constants $T(R_{1})$ and $R_{0}$ do not depend on $\gamma$.
\end{remark}

It follows from Lemma \ref{Absorbing} that the ball $B_{0}=\{x:\left\Vert
x\right\Vert \leq R_{0}\}$ is absorbing for $T_{\gamma}.$ Next, we will obtain
an estimate of the tails of solutions.

\begin{lemma}
\label{Tails}For any $R_{1}>0$, $\varepsilon>0,$ there are $T(\varepsilon
,R_{1})$ and $K(\varepsilon,R_{1})$ such that, any solution $u\left(
\text{\textperiodcentered}\right)  $ of (\ref{EqAut}) starting at $B_{R_{1}}$
satisfies
\[
\int_{\left\vert x\right\vert \geq k}\left\vert u(t,x)\right\vert ^{2}%
dx\leq\varepsilon, \quad\text{ for all }t\geq T,\ k\geq K.
\]

\end{lemma}

\begin{proof}
Let $\theta:\mathbb{R}^{+}\rightarrow\mathbb{R}^{+}$ be a smooth function such
that%
\[
\theta(s)=%
\begin{cases}%
\begin{split}
0,\quad &  \text{ if }~0\leq s\leq\frac{1}{2},\\
1,\quad &  \text{ if }~s\geq1.
\end{split}
\end{cases}
\]
Multiplying the equation in (\ref{EqAut}) by $\theta\left(  \frac{\left\vert
x\right\vert }{k}\right)  u$, we obtain%
\begin{align}
&  \frac{1}{2}\frac{d}{dt}\int_{\mathbb{R}^{m}}\theta\left(  \frac{\left\vert
x\right\vert }{k}\right)  \left\vert u(t,x)\right\vert ^{2}dx+\mu
\int_{\mathbb{R}^{m}}\theta\left(  \frac{\left\vert x\right\vert }{k}\right)
\left\vert u(t,x)\right\vert ^{2}dx\label{EstTail1}\\
&  =-\frac{1}{2}C(m,\gamma)\int_{\mathbb{R}^{m}}\int_{\mathbb{R}^{m}}%
\frac{\left(  \theta\left(  \frac{\left\vert x\right\vert }{k}\right)
u(t,x)-\theta\left(  \frac{\left\vert y\right\vert }{k}\right)  u(t,y)\right)
(u(t,x)-u(t,y))}{\left\vert x-y\right\vert ^{m+2\gamma}}dydx\nonumber\\
&  +\int_{\mathbb{R}^{m}}\theta\left(  \frac{\left\vert x\right\vert }%
{k}\right)  (f(x,u(t,x))+h(x))u(t,x)dx.\nonumber
\end{align}
For the first term of right hand side of (\ref{EstTail1}), by the H\"{o}lder
inequality, we have%
\begin{align*}
&  -\int_{\mathbb{R}^{m}}\int_{\mathbb{R}^{m}}\frac{\left(  \theta\left(
\frac{\left\vert x\right\vert }{k}\right)  u(t,x)-\theta\left(  \frac
{\left\vert y\right\vert }{k}\right)  u(t,y)\right)  (u(t,x)-u(t,y))}%
{\left\vert x-y\right\vert ^{m+2\gamma}}dydx\\
&  \leq-\int_{\mathbb{R}^{m}}\int_{\mathbb{R}^{m}}\frac{\left(  \theta\left(
\frac{\left\vert x\right\vert }{k}\right)  -\theta\left(  \frac{\left\vert
y\right\vert }{k}\right)  \right)  (u(t,x)-u(t,y))u(t,y)}{\left\vert
x-y\right\vert ^{m+2\gamma}}dxdy\\
&  \leq\left\Vert u(t)\right\Vert \left(  \int_{\mathbb{R}^{m}}\left(
\int_{\mathbb{R}^{m}}\frac{\left\vert \theta\left(  \frac{\left\vert
x\right\vert }{k}\right)  -\theta\left(  \frac{\left\vert y\right\vert }%
{k}\right)  \right\vert \left\vert u(t,x)-u(t,y)\right\vert }{\left\vert
x-y\right\vert ^{m+2\gamma}}dx\right)  ^{2}dy\right)  ^{\frac{1}{2}}\\
&  \leq\left\Vert u(t)\right\Vert \left(  \int_{\mathbb{R}^{m}}\int%
_{\mathbb{R}^{m}}\frac{\left\vert \theta\left(  \frac{\left\vert x\right\vert
}{k}\right)  -\theta\left(  \frac{\left\vert y\right\vert }{k}\right)
\right\vert ^{2}}{\left\vert x-y\right\vert ^{m+2\gamma}}dx\int_{\mathbb{R}%
^{m}}\frac{\left\vert u(t,x)-u(t,y)\right\vert ^{2}}{\left\vert x-y\right\vert
^{m+2\gamma}}dxdy\right)  ^{\frac{1}{2}}.
\end{align*}
Since $\left\vert \theta^{\prime}(s)\right\vert \leq L_{1}$ for some $L_{1}%
>0$, by means of the change of variable $z=\frac{x}{k}$, we know that there
exists a constant $L_{2}>0$ such that,
\begin{align*}
&  \int_{\mathbb{R}^{m}}\frac{\left\vert \theta\left(  \frac{\left\vert
x\right\vert }{k}\right)  -\theta\left(  \frac{\left\vert y\right\vert }%
{k}\right)  \right\vert ^{2}}{\left\vert x-y\right\vert ^{m+2\gamma}}dx\\
&  =\frac{1}{k^{2\gamma}}\int_{\mathbb{R}^{m}}\frac{\left\vert \theta\left(
\left\vert z\right\vert \right)  -\theta\left(  \frac{\left\vert y\right\vert
}{k}\right)  \right\vert ^{2}}{\left\vert z-\frac{y}{k}\right\vert
^{m+2\gamma}}dz=\frac{1}{k^{2\gamma}}\int_{\mathbb{R}^{m}}\frac{\left\vert
\theta\left(  \left\vert \xi+\frac{y}{k}\right\vert \right)  -\theta\left(
\frac{\left\vert y\right\vert }{k}\right)  \right\vert ^{2}}{\left\vert
\xi\right\vert ^{m+2\gamma}}d\xi\\
&  \leq\frac{L_{1}^{2}}{k^{2\gamma}}\int_{|\xi|\leq1}\frac{1}{\left\vert
\xi\right\vert ^{m+2\gamma-2}}d\xi+\frac{1}{k^{2\gamma}}\int_{|\xi|\geq1}%
\frac{1}{\left\vert \xi\right\vert ^{m+2\gamma}}d\xi\leq\frac{L_{2}%
}{k^{2\gamma}}.
\end{align*}
Thus, the above estimate and Lemma \ref{Estimates} imply that there exists a
constant $L_{3}>0$, such that
\begin{align}
&  -\frac{1}{2}C(m,\gamma)\int_{\mathbb{R}^{m}}\int_{\mathbb{R}^{m}}%
\frac{\left(  \theta\left(  \frac{\left\vert x\right\vert }{k}\right)
u(t,x)-\theta\left(  \frac{\left\vert y\right\vert }{k}\right)  u(t,y)\right)
(u(t,x)-u(t,y))}{\left\vert x-y\right\vert ^{m+2\gamma}}dydx\nonumber\\[0.8ex]
&  \leq\frac{\sqrt{L_{2}}}{2k^{\gamma}}C(m,\gamma)\left\Vert u(t)\right\Vert
\Vert u(t)\Vert_{\dot{H}^{\gamma}(\mathbb{R}^{m})}\nonumber\\[0.8ex]
&  \leq\frac{L_{3}C(m,\gamma)}{k^{\gamma}}\left(  1+\Vert u(t)\Vert_{\dot
{H}^{\gamma}(\mathbb{R}^{m})}^{2}\right)  . \label{EstTail2}%
\end{align}
For the second term of right hand side of (\ref{EstTail1}), using
(\ref{DissAut}) and the Young inequality, we find
\begin{align}
&  \int_{\mathbb{R}^{m}}\theta\left(  \frac{\left\vert x\right\vert }%
{k}\right)  (f(x,u(t,x))+h(x))u(t,x)dx\label{EstTail3}\\
&  \leq\frac{\mu}{2}\int_{\mathbb{R}^{m}}\theta\left(  \frac{\left\vert
x\right\vert }{k}\right)  \left\vert u(t,x)\right\vert ^{2}dx+\int%
_{\mathbb{R}^{m}}\theta\left(  \frac{\left\vert x\right\vert }{k}\right)
\psi_{1}(x)dx+\frac{1}{2\mu}\int_{\mathbb{R}^{m}}\theta\left(  \frac
{\left\vert x\right\vert }{k}\right)  \left\vert h(x)\right\vert
^{2}dx.\nonumber
\end{align}
Since $\psi_{1}\in L^{1}(\mathbb{R}^{m})$ and $h\in H$, we deduce that for
$\varepsilon^{\prime}>0$, there exists $K(\varepsilon^{\prime})>0$ such that
for all $k\geq K$, we have
\begin{align}
\int_{\mathbb{R}^{m}}\theta\left(  \frac{\left\vert x\right\vert }{k}\right)
\psi_{1}(x)dx  &  \leq\varepsilon^{\prime},\nonumber\\[1ex]
\frac{1}{2\mu}\int_{\mathbb{R}^{m}}\theta\left(  \frac{\left\vert x\right\vert
}{k}\right)  \left\vert h(x)\right\vert ^{2}dx  &  \leq\varepsilon^{\prime
},\nonumber\\[1ex]
\max\left\{  \frac{L_{3}}{k^{\gamma}},\frac{L_{3}C(m,\gamma)}{k^{\gamma}%
}\right\}   &  \leq\varepsilon^{\prime}. \label{EstTail4}%
\end{align}
Collecting (\ref{EstTail1})-(\ref{EstTail4}), it yields%
\[
\frac{d}{dt}\int_{\mathbb{R}^{m}}\theta\left(  \frac{\left\vert x\right\vert
}{k}\right)  \left\vert u(t,x)\right\vert ^{2}dx+\mu\int_{\mathbb{R}^{m}%
}\theta\left(  \frac{\left\vert x\right\vert }{k}\right)  \left\vert
u(t,x)\right\vert ^{2}dx\leq6\varepsilon^{\prime}+2\varepsilon^{\prime
}C(m,\gamma)\Vert u(t)\Vert_{\dot{H}^{\gamma}(\mathbb{R}^{m})}^{2}.
\]

Multiplying the above inequality by $e^{\mu s}$ and integrating it over
$(0,t)$, together with estimate (\ref{EstsolGamma}), we infer there exists a
positive constant $L_{4}$ such that,
\begin{align*}
\int_{\mathbb{R}^{m}}\theta\left(  \frac{\left\vert x\right\vert }{k}\right)
\left\vert u(t,x)\right\vert ^{2}dx  &  \leq e^{-\mu t}\int_{\mathbb{R}^{m}%
}\theta\left(  \frac{\left\vert x\right\vert }{k}\right)  \left\vert
u_{0}(x)\right\vert ^{2}dx\\
&  ~~+\frac{6\varepsilon^{\prime}}{\mu}+2\varepsilon^{\prime}C(m,\gamma
)\int_{0}^{t}e^{-\mu(t-s)}\Vert u(s)\Vert_{\dot{H}^{\gamma}(\mathbb{R}^{m}%
)}^{2}ds\\
&  \leq e^{-\mu t}\int_{\mathbb{R}^{m}}\theta\left(  \frac{\left\vert
x\right\vert }{k}\right)  \left\vert u_{0}(x)\right\vert ^{2}dx+L_{4}%
\varepsilon^{\prime}.
\end{align*}
We can choose $T(\varepsilon^{\prime})>0$ such that $e^{-\mu t}\int%
_{\mathbb{R}^{m}}\theta\left(  \frac{\left\vert x\right\vert }{k}\right)
\left\vert u_{0}(x)\right\vert ^{2}dx\leq\varepsilon^{\prime}$ for all $t\geq
T$, and $\varepsilon^{\prime}$ satisfies $\left(  1+L_{4}\right)
\varepsilon^{\prime}<\varepsilon$. Therefore,
\[
\int_{\left\vert x\right\vert \geq k}\left\vert u(t,x)\right\vert ^{2}%
dx\leq\int_{\mathbb{R}^{m}}\theta\left(  \frac{\left\vert x\right\vert }%
{k}\right)  \left\vert u(t,x)\right\vert ^{2}dx\leq\varepsilon,\quad\forall
t\geq T, \forall k\geq K.
\]
The proof is complete.
\end{proof}

\begin{lemma}
\label{AC} The operator $T_{\gamma}$ is asymptotically compact, that is, for
any bounded set $B\in H$, every sequence $y_{\gamma}^{n}=T_{\gamma}%
(t_{n},x_{n})$ with $x_{n}\in B$ is relatively compact in $H$ when
$t_{n}\rightarrow+\infty$.
\end{lemma}

\begin{proof}
In view of (\ref{EstsolGamma}), the sequence $\{y_{\gamma}^{n}\}$ is bounded
in $H$. Hence, up to a subsequence (relabeled the same), $y_{\gamma}%
^{n}\rightarrow y$ weakly in $H.$ We need to prove that the convergence is in
fact strong. As $y_{\gamma}^{n}=T_{\gamma}(1,T_{\gamma}(t_{n}-1,x_{n}))$,
there exist solutions $u_{\gamma}^{n}\left(  \text{\textperiodcentered
}\right)  $ such that $u_{\gamma}^{n}(0)=z_{\gamma}^{n}=T_{\gamma}%
(t_{n}-1,x_{n})$ and $u_{\gamma}^{n}(1)=y_{\gamma}^{n}$. The sequence
$z_{\gamma}^{n}$ converges weakly in $H$ to some $z$ as well.

It follows from (\ref{EstsolGamma})-(\ref{EstSol2}) that there is $M>0$ such
that%
\[
\sup_{t\in\lbrack0,1]}\ \left\Vert u_{\gamma}^{n}(t)\right\Vert \leq M,
\]%
\[
\int_{0}^{1}\left(  C(m,\gamma)\Vert u_{\gamma}^{n}(s)\Vert_{\dot{H}^{\gamma
}(\mathbb{R}^{m})}^{2}+2\beta\left\Vert u_{\gamma}^{n}(s)\right\Vert _{p}%
^{p}\right)  ds\leq M.
\]
Thus, $\{u_{\gamma}^{n}\}$ is bounded in $L^{\infty}(0,1;H)\cap L^{p}%
(0,1;L^{p}(\mathbb{R}^{m}))\cap L^{2}(0,1;V_{\gamma})$. It also follows that
$\{A^{\gamma}(u_{\gamma}^{n})\}$ is bounded in $L^{2}(0,1;V_{\gamma}^{\ast})$.
Moreover, by (\ref{GrowthAut}), we deduce that $\{F(u_{\gamma}^{n})\}$ is
bounded in $L^{q}(0,1;L^{q}(\mathbb{R}^{m}))$. Hence, $\left\{  \dfrac
{du_{\gamma}^{n}}{dt}\right\}  $ is bounded in $L^{q}(0,1;L^{q}(\mathbb{R}%
^{m}))+L^{2}(0,1;V_{\gamma}^{\ast})$. Then, using the Aubin-Lions lemma, there
exists a function $u_{\gamma}\in L^{\infty}(0,1;H)\cap L^{p}(0,1;L^{p}%
(\mathbb{R}^{m}))\cap L^{2}(0,1;V^{\ast})$ with $\dfrac{du_{\gamma}}{dt}\in
L^{q}(0,1;L^{q}(\mathbb{R}^{m}))+L^{2}(0,1;V_{\gamma}^{\ast})$ and $\chi\in
L^{q}(0,1;L^{q}(\mathbb{R}^{m}))$, such that, up to a subsequence (relabeled
the same), the following convergences hold:%
\begin{equation}
u_{\gamma}^{n}\rightarrow u_{\gamma}\text{ weak-star in }L^{\infty}(0,1;H),
\label{Conv1}%
\end{equation}%
\begin{equation}
u_{\gamma}^{n}\rightarrow u_{\gamma}\text{ weakly in }L^{2}(0,1;V_{\gamma}),
\label{Conv2}%
\end{equation}%
\begin{equation}
u_{\gamma}^{n}\rightarrow u_{\gamma}\text{ weakly in }L^{p}(0,1;L^{p}%
(\mathbb{R}^{m})), \label{Conv3}%
\end{equation}%
\begin{equation}
A^{\gamma}(u_{\gamma}^{n})\rightarrow A^{\gamma}(u_{\gamma})\text{ weakly in
}L^{2}(0,1;V_{\gamma}^{\ast}), \label{Conv4}%
\end{equation}%
\begin{equation}
\frac{du_{\gamma}^{n}}{dt}\rightarrow\frac{du_{\gamma}}{dt}\text{ weakly in
}L^{q}(0,1;L^{q}(\mathbb{R}^{m}))+L^{2}(0,1;V_{\gamma}^{\ast}), \label{Conv5}%
\end{equation}%
\begin{equation}
F(u_{\gamma}^{n})\rightarrow\chi\text{ weakly in }L^{q}(0,1;L^{q}%
(\mathbb{R}^{m})). \label{Conv6}%
\end{equation}

We need to prove that $u_{\gamma}\left(  \text{\textperiodcentered}\right)  $
is a weak solution of problem (\ref{EqAut}) satisfying $u_{\gamma}\left(
0\right)  =z$ and $u_{\gamma}(1)=y$. To this end, for $R>0$, we define the
projections $u_{\gamma,R}^{n}=L_{R}u_{\gamma}^{n}$, where $L_{R}:H\rightarrow
L^{2}(B_{R})$ is defined by $L_{R}v(x)=v(x)$ for a.a. $x\in B_{R}$. It is easy
to see that the above convergences (\ref{Conv1})-(\ref{Conv6}) still hold, if
we replace $u_{\gamma}^{n}$ by $u_{\gamma,R}^{n},\ u_{\gamma}$ by
$L_{R}u_{\gamma}$, $\mathbb{R}^{m}$ by $B_{R}$, $H$ by $L^{2}(B_{R})$,
$V_{\gamma}$ by $H^{\gamma}(B_{R})$, $V_{\gamma}^{\ast}$ by $\left(
H^{\gamma}(B_{R})\right)  ^{\ast}$ and $\chi$ by $L_{R}\chi$, respectively.

Since the embedding $H_{\gamma}(B_{R})\subset L^{2}(B_{R})$ is compact (cf.
\cite{Nezza}) and the embedding $L^{2}(B_{R})\subset(H^{\gamma}(B_{R})\cap
L^{p}(B_{R}))^{\ast}$ is continuous, a standard Compactness Theorem
\cite[Theorem 8.1]{Robinson} implies that%
\begin{equation}
u_{\gamma,R}^{n}\rightarrow L_{R}u_{\gamma}\text{ strongly in }L^{2}%
(0,1;L^{2}(B_{R})), \label{Converg7}%
\end{equation}%
\begin{equation}
u_{\gamma,R}^{n}(t,x)\rightarrow L_{R}u_{\gamma}(t,x)~\text{ for a.a. }%
t\in\left(  0,1\right)  , x\in B_{R}. \label{Converg8}%
\end{equation}
Hence, $f(x,u_{\gamma,R}^{n}(t,x))\rightarrow f(x,L_{R}u_{\gamma}(t,x))$ for
a.a. $t\in(0,1)$ and $x\in B_{R}.$ \cite[ Lemma 1.3]{Lions} implies that
$L_{R}\chi=F_{R}(L_{R}u_{\gamma})$, where $F_{R}:L^{p}(B_{R})\rightarrow
L^{q}(B_{R})$ is the corresponding Nemytskii operator associated to $f$. In a
standard way (see e.g., \cite[Lemma 8]{MorVal}), one can show that $u_{\gamma
}$ satisfies the equality in (\ref{EqAut}) in the sense of distributions. Then
equality (\ref{EqualitySol}) is true as well, so $u_{\gamma}$ is a solution of
(\ref{EqAut}).

Further, we need to prove that $u_{\gamma}\left(  0\right)  =z$ and
$u_{\gamma}(1)=y$. As the embedding $L^{2}(B_{R})\subset(H^{\gamma}(B_{R})\cap
L^{p}(B_{R}))^{\ast}$ is compact and $\dfrac{du_{\gamma,R}^{n}}{dt}$ is
bounded in $L^{q}(0,1;(H^{\gamma}(B_{R})\cap L^{p}(B_{R}))^{\ast})$, the
Ascoli-Arzel\`{a} Theorem implies that $u_{\gamma,R}^{n}\rightarrow
L_{R}u_{\gamma}$ strongly in $C([0,1],(H^{\gamma}(B_{R})\cap L^{p}%
(B_{R}))^{\ast})$. From this and (\ref{EstsolGamma}), we deduce immediately
that $u_{\gamma,R}^{n}(0)\rightarrow L_{R}u_{\gamma}(0),\ u_{\gamma,R}%
^{n}(1)\rightarrow L_{R}u_{\gamma}(1)$ weakly in $L^{2}(B_{R})$. Hence,
$L_{R}z=L_{R}u_{\gamma}(0)$ and $L_{R}y=L_{R}u_{\gamma}(1)$ for any $R>0$.
Thus, the result follows.

On the one hand, integrating in (\ref{EstSol2}) over $(0,t)$ with $t>0$, we
find that there exists a constant $C>0$ such that the functions
$J(t)=\left\Vert u_{\gamma}(t)\right\Vert ^{2}-Ct,\ J_{n}(t)=\left\Vert
u_{\gamma}^{n}(t)\right\Vert ^{2}-Ct$ are non-increasing. On the other hand,
in terms of (\ref{Converg7}), for any $R>0$, we have that $L_{R}u_{\gamma}%
^{n}(t)\rightarrow L_{R}u_{\gamma}(t)$ in $L^{2}(B_{R})$ for a.a. $t\in\left(
0,1\right)  $. For any $\varepsilon>0$, Lemma \ref{Tails} implies the
existence of $R_{2}(\varepsilon)>0$ and $N_{1}(\varepsilon)>0$ such that%
\[
\int_{\left\vert x\right\vert >R_{2}}\left\vert u_{\gamma}^{n}(s,x)\right\vert
^{2}dx\leq\varepsilon,\quad\text{ for all }s\in\lbrack0,1]\text{ and }n\geq
N_{1}.
\]
Take a sequence $t_{m}\rightarrow1^{-}$ such that $L_{R}u_{\gamma}^{n}%
(t_{m})\rightarrow L_{R}u_{\gamma}(t_{m})$ in $L^{2}(B_{R})$ for any $t_{m}%
\in(0,1)$. Then by the monotonicity and the continuity of $J_{n}$ and $J$, we
infer that there is $N_{2}(\varepsilon)\geq N_{1}(\varepsilon)$ such that%
\begin{align*}
J_{n}(1)-J(1)  &  =J_{n}(1)-J_{n}(t_{m})+J_{n}(t_{m})-J(t_{m})+J(t_{m})-J(1)\\
&  \leq2\varepsilon+\left\vert \int_{\left\vert x\right\vert \leq R_{2}%
}\left\vert u_{\gamma}^{n}(t_{m},x)\right\vert ^{2}dx-\int_{\left\vert
x\right\vert \leq R_{2}}\left\vert u_{\gamma}(t_{m},x)\right\vert
^{2}dx\right\vert \\
&  +\int_{\left\vert x\right\vert >R_{2}}\left\vert u_{\gamma}^{n}%
(t_{m},x)\right\vert ^{2}dx+\int_{\left\vert x\right\vert >R_{2}}\left\vert
u_{\gamma}(t_{m},x)\right\vert ^{2}dx\\
&  \leq5\varepsilon.
\end{align*}
Hence, $\limsup_{n\rightarrow\infty}J_{n}(1)\leq J(1)$, which implies that
$\limsup_{n\rightarrow\infty}\left\Vert u_{\gamma}^{n}(1)\right\Vert
\leq\left\Vert u_{\gamma}(1)\right\Vert $. Since $\left\Vert y\right\Vert
\leq\liminf_{n\rightarrow\infty}\left\Vert y_{n}\right\Vert $, we obtain that
$y_{n}\rightarrow y$ strongly in $H.$
\end{proof}

\begin{theorem}
$T_{\gamma}$ has a connected global attractor $\mathcal{A}_{\gamma}$, which is
characterized by%
\[
\mathcal{A}_{\gamma}=\{\phi\left(  0\right)  :\phi\text{ is a bounded complete
trajectory of }T_{\gamma}\}. \label{Charac}%
\]
Moreover, it is the minimal closed set attracting every bounded set.
\end{theorem}

\begin{proof}
In view of Lemmas \ref{Absorbing} and \ref{AC}, it follows from standard
results (see e.g. \cite[Theorem 3.1]{Lad}).
\end{proof}

\begin{lemma}
\label{AttrBounded}The set $\cup_{\gamma\in\left(  0,1\right)  }%
\mathcal{A}_{\gamma}$ is bounded in $H$.
\end{lemma}

\begin{proof}
It follows from the inclusion $\mathcal{A}_{\gamma}\subset B_{0}$ for any
$\gamma\in\left(  0,1\right)  .$
\end{proof}

\section*{Acknowledgements}

This research was supported by FEDER and the Spanish Ministerio de Ciencia e
Innovaci\'{o}n under projects PDI2021-122991-NB-C21 and PGC2018-096540-B-I00,
Junta de Andaluc\'{\i}a (Spain) under project P18-FR-4509 and by the
Generalitat Valenciana, project PROMETEO/2021/063, and the Nature Science
Foundation of Jiangsu Province (Grant No. BK20220233).

\section*{Data availability}
No data has been used in the development of the research in this paper.

\end{document}